\documentclass[preprint]{elsarticle}
\usepackage{amsmath,amssymb,mathrsfs}
\usepackage{graphicx,color,url}
\usepackage{amsthm}
\usepackage{subfig}
\usepackage[margin=1in]{geometry}
\usepackage[linesnumbered,ruled,vlined]{algorithm2e}
\usepackage{layouts}

\newcommand{\fundef}[3]{#1:#2\to #3}

\newcommand{\xb}{\mathbf{x}}
\newcommand{\ub}{\mathbf{u}}

\newcommand{\R}{\mathbb{R}}
\newcommand{\Q}{\mathbb{Q}}

\newcommand{\bz}{B\'{e}zier }

\newcommand{\po}{\partial\Omega}
\newcommand{\dd}{\mathrm{d}\Omega}

\newcommand{\Tb}{\mathbf{T}}
\newcommand{\Tc}{\mathcal{T}}
\newcommand{\Qb}{\mathbf{Q}}
\newcommand{\Qc}{\mathcal{Q}}
\newcommand{\jacob}[1]{\det\left(\nabla#1\right)}

\theoremstyle{definition}
\newtheorem{definition}{Definition}

\newcounter{myremark}
\newenvironment{myremark}[1][]{\refstepcounter{myremark}\par\medskip
   \noindent \textbf{Remark~\themyremark.
	 #1} \rmfamily}{\medskip}

\SetCommentSty{mycommfont}
\SetKwInput{KwInput}{Input}                
\SetKwInput{KwOutput}{Output}              

\biboptions{sort&compress}

\begin{document}

\sloppy

\begin{frontmatter}

\journal{Computer Methods in Applied Mechanics and Engineering}

\title{Robust Numerical Integration on Curved Polyhedra \\Based on Folded Decompositions}

\author[1]{Pablo Antolin\corref{cor1}}
\author[1]{Xiaodong Wei}
\author[1,2]{Annalisa Buffa}
\cortext[cor1]{Corresponding author, pablo.antolin@epfl.ch}
\address[1]{Institute of Mathematics, \'Ecole Polytechnique F\'ed\'erale de Lausanne, 1015 Lausanne, Switzerland}
\address[2]{Istituto di Matematica Applicata e Tecnologie Informatiche ``Enrico Magenes'' del CNR, Pavia, Italy}

\begin{abstract}
We present a novel method to perform numerical integration over curved polyhedra enclosed by high-order parametric surfaces.
Such a polyhedron is first decomposed into a set of triangular and/or rectangular pyramids, whose certain faces correspond to the given parametric surfaces.
Each pyramid serves as an integration cell with a geometric mapping from a standard parent domain (e.g.,  a unit cube), where the tensor-product Gauss quadrature is adopted.
As no constraint is imposed on the decomposition, certain resulting pyramids may intersect with themselves, and thus their geometric mappings may present negative Jacobian values.
We call such cells the \emph{folded cells} and refer to the corresponding decomposition as a folded decomposition.
We show that folded cells do not cause any issues in practice as they are only used to numerically compute certain integrals of interest.
The same idea can be applied to planar curved polygons as well.
We demonstrate both theoretically and numerically that folded cells can retain the same accuracy as the cells with strictly positive Jacobians.
On the other hand, folded cells allow for a much easier and much more flexible decomposition for general curved polyhedra, on which one can robustly compute integrals.
In the end, we show that folded cells can flexibly and robustly accommodate real-world complex geometries by presenting several examples in the context of immersed isogeometric analysis, where involved sharp features can be well respected in generating integration cells.
\end{abstract}

%
%
%


\end{frontmatter}

\section{Introduction}


Accurate, efficient, and robust numerical integration methods over curved polygons/polyhedra are required in various applications across computer-aided engineering (CAE), computer-aided design (CAD) \cite{ref:krishnamurthy11}, and computer graphics and animation \cite{ref:guendelman03}.
Particularly in CAE, novel integration schemes are often needed when the mesh for approximating solutions does not conform to the geometric representation, such as extended finite element methods \cite{ref:moes99}, cut finite element methods (cutFEM) \cite{ref:burman15}, virtual element methods \cite{ref:veiga13a}, finite cell methods (FCM) \cite{ref:parvizian07}, discontinuous-Galerkin-based immersed boundary methods \cite{ref:rangarajan09}, and immersed isogeometric analysis (IGA) \cite{ref:kamensky15}.
The key idea of such non-conformal methods is to decouple geometric representation from the discretization for solution approximation to significantly simplify mesh generation.
For instances, in FCM and immersed IGA, it involves embedding the geometric representation into a simple background grid (e.g., a Cartesian grid), leading to cut (or trimmed) elements of complex shapes.
Properly dealing with such cut elements poses several challenges such as numerical integration, stabilization regardless of cutting \cite{ref:marussig17, ref:johansson19, ref:puppi20, ref:wei21}, imposition of the Dirichlet boundary conditions \cite{ref:ruess13}, and so forth.
We focus on the numerical integration in this paper.


There exist several major groups of methods to perform numerical integration over curved polygons/polyhedra, such as (1) application of the divergence theorem to reduce the dimension of integrals, (2) moment fitting, and (3) domain decomposition into simplices.
First, the divergence theorem can be applied to convert an integral to its low-dimension counterpart \cite{ref:timmer80, ref:sudhakar14, ref:chin20, ref:gunderman20, ref:antolin21}.
For instance, a volume integral can be reduced to a surface integral, and even to a point evaluation \cite{ref:chin20, ref:antolin21}.
To be able to apply the divergence theorem, it is critical to find a suitable antiderivative for the integrand.
Accordingly, one needs to either assume integrands to be in a certain form (e.g., homogeneous functions) to have analytic antiderivatives \cite{ref:chin20} or employ certain numerical methods to obtain approximations \cite{ref:gunderman20}.

Second, moment fitting solves a system of equations to find a quadrature rule (i.e., the locations and/or weights of quadrature points) that can exactly integrate certain polynomial basis functions \cite{ref:xiao10, ref:mousavi11, ref:sudhakar13, ref:muller13, ref:joulaian16}.
When point locations are taken as a priori, the system becomes linear and the unknowns are only weights.
Recently, moment fitting was further explored in nonlinear problems \cite{ref:hubrich2017, ref:hubrich19, ref:bui20}.
The advantage of moment fitting is to retain the minimum number of integration points at the cost of solving a linear system for each cut element.
In fact, it needs to be used together with other groups of methods to accurately compute the right-hand-side vectors.
In other words, moment fitting moves the computation of integrals over complex domains to the preprocessing step such that later it can be used in a fast way.


Third, decomposition-based methods represent a popular way of computing integrals over complex domains.
The basic idea is to partition a given domain into manageable ``cells" where standard quadrature can be applied.
Such methods can be divided into low-order and high-order ones.
A typical example of low-order decomposition is adaptive quad/octree subdivision widely adopted in FCM \cite{ref:schillinger15}, where a cut element is adaptively subdivided into sub-cells to capture the domain boundary for integration.
Essentially, it is merely a piecewise constant approximation to the given geometry.
Such a low-order approximation leads to a large number of quadrature points aggregated around the boundary and thus renders it computationally expensive especially in nonlinear problems or in 3D \cite{ref:duster08,ref:divi2020, ref:peto20}.
On the other hand, high-order decomposition features a high-order approximation to the geometry of interest, and thus it yields much fewer quadrature cells than the adaptive subdivision; see \cite{ref:marco15, ref:kudela16, ref:im18, ref:antolin19a, ref:massarwi19} for several methods that work in 3D.
Due to its high accuracy per quadrature point, it is particularly suitable for high-order methods such as immersed IGA.
However, a high-order decomposition generally involves a high-order triangulation \cite{ref:antolin19a} or other sophisticated partitions \cite{ref:marco15}.
The procedure can be complicated and time-consuming especially when the resulting integration cells are required to have positive Jacobian.
However, we will show that the positive-Jacobian constraint is in fact not necessary for the integration purpose.

In this paper, we propose a novel decomposition-based method that allows also for integration cells that present negative Jacobian values. We call such cells \emph{folded cells} and refer to the corresponding decomposition as a \emph{folded decomposition}.
In other words, no constraint is imposed to the Jacobians of the integration cells when we create a decomposition for a given domain.
Allowing for folded cells greatly increases the flexibility and guarantees that the decomposition is always successful.
The method works for both low-order and high-order decompositions.
We show that the performance of folded decompositions is comparable to that of strictly positive-Jacobian cells in terms of integration accuracy and convergence.
Moreover, when applied to solving linear elliptic partial differential equations, our strategy yields expected error convergence rates, and due to its intrinsic flexibility, it can robustly accommodate real-world complex geometries.
We note that another independent work using folded cells was recently reported in \cite{ref:chin20a}, which, however, is limited to 2D domains and did not theoretically show how such cells work.

The reminder of the paper is organized as follows.
Section \ref{sec:ps} presents the problem setting.
Next in Sections \ref{sec:nj2d} and \ref{sec:nj3d}, we discuss the folded decomposition for both curved polygons and curved polyhedra.
Section \ref{sec:iga} discusses how to apply the proposed method to immersed IGA.
Numerical examples are then presented in Sections \ref{sec:int_poly} and \ref{sec:pde}, where Section \ref{sec:int_poly} is devoted to demonstrating the efficacy and robustness of folded decomposition, whereas Section \ref{sec:pde} shows that folded decomposition can robustly accommodate real-world complex geometries in engineering simulations.
We conclude our work in Section \ref{sec:con}.

\section{Problem setting: Integration over general domains}
\label{sec:ps}

We aim to find an accurate and robust quadrature rule to approximate the integral of a certain function $f(\xb)$ over a compact domain $\Omega\subset\R^d$ ($d=2,3$, the spatial dimension),
\begin{equation}
\int_{\Omega} f(\xb) \, \dd.
\label{eq:integral}
\end{equation}
In finite element analysis (FEA), $f(\xb)$ usually represents quantities contributing to mass or stiffness matrices.
The challenge of computing \eqref{eq:integral} originates from the complexity of the domain $\Omega$.

In this paper, $\Omega$ is assumed to be a curved polygonal domain in 2D or a curved polyhedral domain in 3D.
More precisely, a curved polygon or polyhedron is a domain $\Omega$ ($d=2$ for a plane or $d=3$ a volume) such that $\po$ is characterized by the union of a set of parameterized edges (for $d=2$), $\{\mathbf{C}_i(u)\in\R^2,\, u\in [0,1],\, i=1,\ldots,N_e\}$, or a set of parameterized faces (for $d=3$), $\{\mathbf{S}_i(u,v)\in\R^3,\, (u,v)\in\mathcal{D}_i\subseteq [0,1]^2,\,i=1,\ldots,N_f\}$. Note that $\mathcal{D}_i$, the definition domain of $\mathbf{S}_i(u,v)$, generally is a subset of the entire parametric domain $[0,1]^2$. 

The aim of this paper is to design robust and cheap approaches to compute Eq. \eqref{eq:integral} for a known regular function $\fundef{f}{\R^d}{\R}$. For the sake of clarity, let us introduce a few terms we use in this paper.

\begin{definition}[Cell]
We call a compact set $\Tc\subset \R^d$ a \emph{cell} if it is defined by a mapping $\fundef{\Tb}{[0,1]^d}{\Tc}$, i.e., $\xb=\Tb(\ub)$, such that its Jacobian, $\mathrm{Jac}_{\Tc}(\xb) \equiv\jacob{\Tb(\ub)}$, is well defined everywhere. The interior and the boundary of $\Tc$ are denoted as $\mathring{\Tc}$ and $\partial\Tc$, respectively, and clearly $\Tc=\mathring{\Tc}\cup\partial\Tc$.
\end{definition}

\begin{definition}[$J^+$ cell]
A cell $\Tc$ is called a \emph{$J^+$ cell} if it has $\mathrm{Jac}_{\Tc}\geq0$ everywhere and $\max(\mathrm{Jac}_{\Tc})>0$.
\end{definition}


\begin{definition}[$J^+$ decomposition]
A set of non-overlapping cells, $\{\Tc_\ell\}_{\ell=1}^L$ where $\mathring{\Tc_\ell} \cap \mathring{\Tc_{\ell'}} =\varnothing$ for any $\ell\neq\ell'$, is called a \emph{$J^+$ decomposition} of $\Omega$ if $\Omega=\cup_{\ell=1}^L \Tc_\ell$ and every $\Tc_\ell$ is a $J^+$ cell.
\end{definition}

With a $J^+$ decomposition, Eq. \eqref{eq:integral} can be computed as
\begin{equation}
\int_{\Omega} f(\xb) = \sum_{\ell=1}^L \int_{\Tc_\ell} f(\xb)
= \sum_{\ell=1}^L \int_{[0,1]^d} f(\Tb_\ell(\ub)) \jacob{\Tb_\ell (\ub)}
\approx \sum_{\ell=1}^L \sum_{q=1}^{N_g} f(\Tb_\ell (\ub_q^g)) \jacob{\Tb_\ell (\ub_q^g)} \, w_{q}^g,
\label{eq:decompquad}
\end{equation}
where $N_g$ is the number of quadrature points, $\ub_q^g$ and $w_q^g$ are tensor-product Gauss quadrature points and weights, respectively, defined in $[0,1]^d$. 
Letting
\begin{equation}
\xb_q^\ell := \Tb_\ell (\ub_q^g) \quad \text{and} \quad  w_q^\ell := \det \left( \nabla \Tb_\ell (\ub_q^g) \right) w_{q}^g,
\end{equation}
we obtain the typical quadrature formula for a decomposition-based method,
\begin{equation}
\int_{\Omega} f(\xb) \approx \sum_{\ell=1}^L \sum_{q=1}^{N_g} f(\xb_q^\ell) \, w_q^\ell,
\label{eq:gptw}
\end{equation}
where we notice that all the weights $w_q^\ell$ are strictly positive.

In this paper, we aim at generalizing this formula to accommodate more general decompositions by allowing for ``folded" cells, i.e., cells that may present negative Jacobian values. Equivalently speaking, the weights $w_q^\ell$ in Eq. \eqref{eq:gptw} may be negative. More specifically, we introduce the following definitions to facilitate our discussion.

\begin{definition}[$J^-$ cell]
A cell $\Tc$ is called a \emph{$J^-$ cell} if it has $\mathrm{Jac}_{\Tc}\leq 0$ everywhere and $\min(\mathrm{Jac}_{\Tc})<0$.
\label{def:j-cell}
\end{definition}

\begin{definition}[Folded cell]
We say that a cell $\Tc$ is a \emph{folded cell} if it has $\min(\mathrm{Jac}_{\Tc})<0$ in $\mathring{\Tc}$, but $\Tc$ can be split into a finite number of subsets, $\Tc = \cup_{i=1}^{N} \Qc_{i}$, such that every $\Qc_{i}$ is either a $J^+$ cell or a $J^-$ cell. 
\label{def:foldcell}
\end{definition}

\begin{figure}[htb]
\centering
\includegraphics[width=0.8\textwidth]{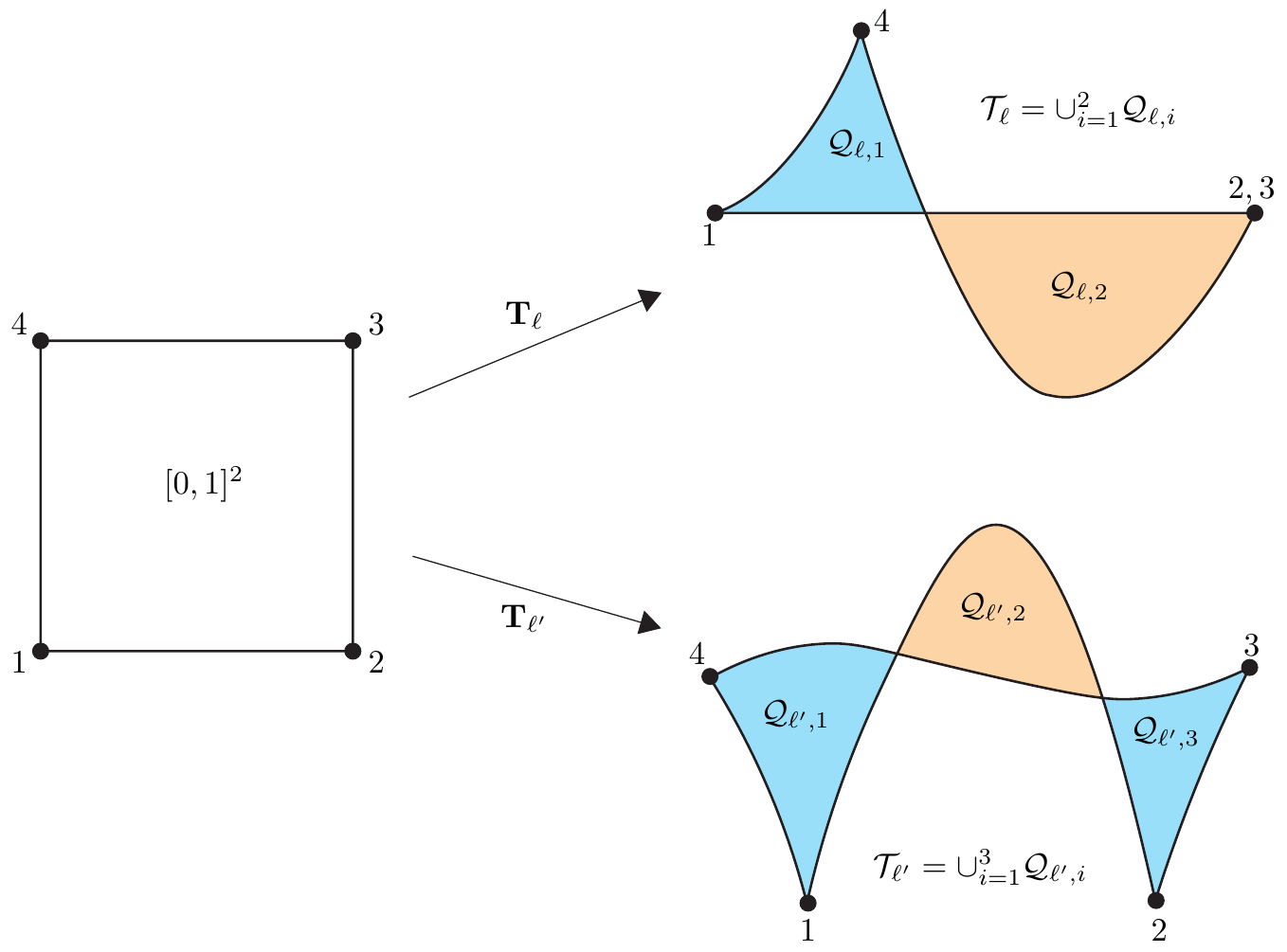}
\caption{Split of folded cells into subsets, inside which the Jacobian does not change sign. The blue and orange regions represent $J^+$ and $J^-$ cells, respectively.}
\label{fig:negative_decomp}
\end{figure}

According to Definitions \ref{def:j-cell} and \ref{def:foldcell}, we note that $J^-$ cells are a special case of folded cells.
A few examples of folded cells are shown in Fig. \ref{fig:negative_decomp}.

\begin{definition}[Folded decomposition]
A set of overlapping cells, $\{\Tc_\ell\}_{\ell=1}^L$ where there exists $\ell\neq\ell'$ such that $\mathring{\Tc_\ell}\cap\mathring{\Tc_{\ell'}}\neq\varnothing$, is called a \emph{folded decomposition} of $\Omega$ if $\Omega=\overline{\sqcup_{\ell=1}^L \mathring{\Tc_\ell}}$ and there exists at least one cell that is a folded cell, where the upper bar denotes the closure of a domain,
\begin{equation}
\sqcup_{\ell=1}^L \mathring{\Tc_\ell} := \{\xb \in \cup_{\ell=1}^L \mathring{\Tc_\ell}: \quad \sum_{\ell=1}^L \mathrm{Count}_{\Tc_\ell}(\xb)=1 \},
\label{eq:signunion}
\end{equation}
and 
\begin{equation}
\mathrm{Count}_{\Tc_\ell}(\xb):=
\begin{cases}
1 & \text{if } \xb\in \mathring{\Tc_\ell} \text{ and } \mathrm{Jac}_{\Tc_\ell}(\xb) > 0,\\
-1 & \text{if } \xb\in \mathring{\Tc_\ell} \text{ and } \mathrm{Jac}_{\Tc_\ell}(\xb) < 0,\\
0 & \text{otherwise.}
\end{cases}
\label{eq:sign}
\end{equation}
\label{def:folddecomp}
\end{definition}

Eqs. (\ref{eq:signunion}, \ref{eq:sign}) in Definition \ref{def:folddecomp} essentially collect the points that only constitute the regions of ``net" positive Jacobians. For example, if a point $\xb$ only belongs to a $J^+$ cell, then $\sum_{\ell=1}^L \mathrm{Count}_{\Tc_\ell}(\xb)=1$ holds and thus $\xb \in \sqcup_{\ell=1}^L \mathring{\Tc_\ell}$. On the other hand, if a point $\xb$ belongs to two cells, $\Tc_\ell$ and $\Tc_{\ell'}$, and moreover $\mathrm{Jac}_{\Tc_\ell}(\xb)>0$ and $\mathrm{Jac}_{\Tc_{\ell'}}(\xb)<0$, then $\sum_{\ell=1}^L \mathrm{Count}_{\Tc_\ell}(\xb)=0$ and thus $\xb \not\in \sqcup_{\ell=1}^L \mathring{\Tc_\ell}$. In other words, the contributions from the positive- and negative-Jacobian regions are canceled out at this same point $\xb$. We will discuss this mechanism in detail in the next section.

In what follows, we provide an efficient algorithm to compute possibly folded decomposition that never fails. 
We start our discussion for $d=2$ in Section \ref{sec:nj2d} and then we move to $d=3$ in Section \ref{sec:nj3d}.

\section{Folded decompositions for $d=2$}
\label{sec:nj2d}


In this section, we introduce two approaches to split 2D domains to generate folded decompositions.

\subsection{Quadrilateral-dominant decomposition}
\label{sec:quad_domi}

\begin{figure}[htb]
\centering
\begin{tabular}{cc}
\includegraphics[width=0.35\textwidth]{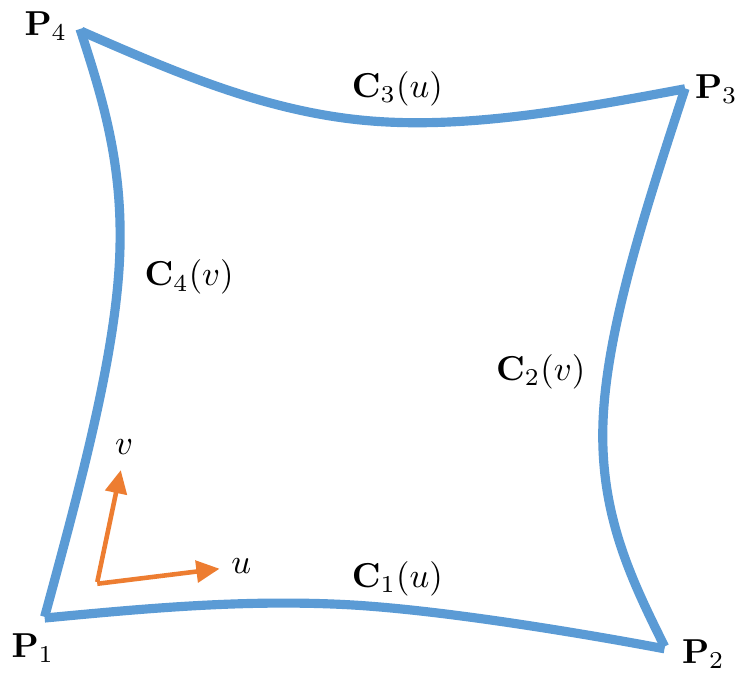} & \hspace{+2em}
\includegraphics[width=0.35\textwidth]{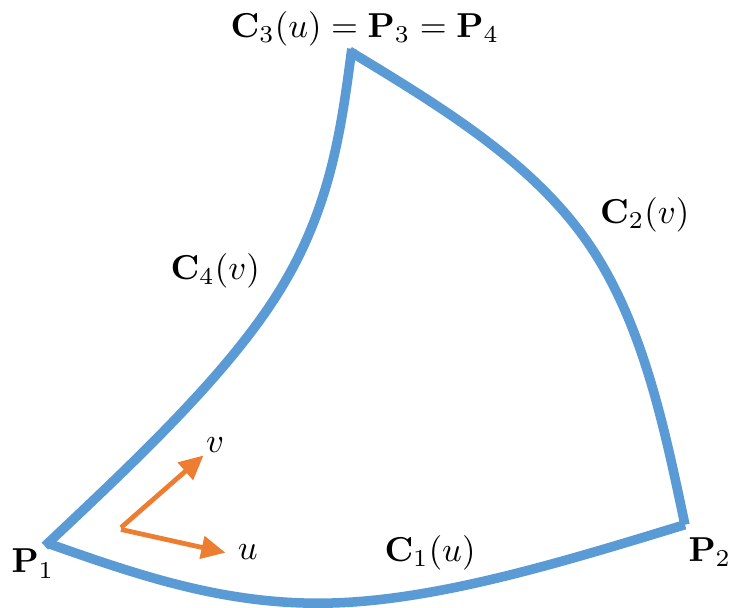} \\
(a) & \hspace{+2em} (b)\\
\end{tabular}
\caption{Boundary curves and corners of Coons patches, where $\mathbf{C}_3(u)$ in (b) is degenerated.}
\label{fig:coons}
\end{figure}

The quadrilateral-dominant decomposition recursively splits a given polygon into quadrilaterals and/or triangles, which provides a means to minimize the number of resulting cells.
Clearly, no decomposition is needed when a polygon is a quadrilateral or a triangle.
When a polygon is a quadrilateral, a Coons patch $\mathbf{T}(u,v)$ \cite{ref:coons} is created to interpolate its 4 edges,
\begin{equation}
\begin{aligned}
\mathbf{T}(u,v) = & (1-v) \mathbf{C}_1(u) + v \mathbf{C}_3(u) + (1-u) \mathbf{C}_4(v) + u \mathbf{C}_2(v) \\
& - (1-u)(1-v) \mathbf{P}_1 - u(1-v) \mathbf{P}_2 - u v \mathbf{P}_3 - (1-u)v \mathbf{P}_4,  \\
\end{aligned}
\label{eq:cell2d}
\end{equation}
where $(u,v)\in [0,1]^2$, and $\fundef{\mathbf{C}_i}{[0,1]}{\R^2}$ and $\mathbf{P}_i\in\R^2$ ($i=1,\ldots,4$) are the corresponding curves and corners, respectively; see Fig.\ \ref{fig:coons}(a).
When a polygon is a triangle, a Coons patch can still be constructed with a degenerated curve (i.e., a zero-length curve); see Fig.\ \ref{fig:coons}(b).
Such (degenerate) Coons patches serve as cells for integration.

\begin{figure}[htb]
\centering
\begin{tabular}{ccc}
\includegraphics[width=0.25\textwidth]{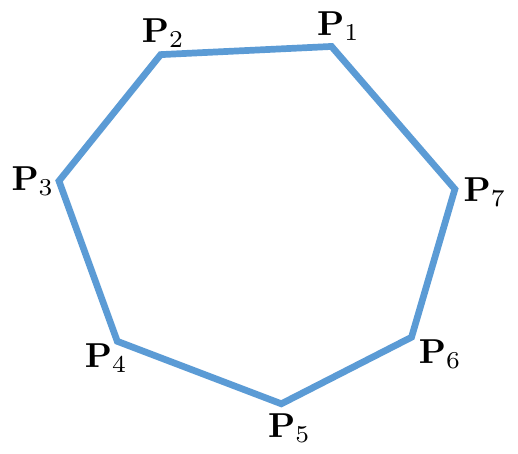} & 
\includegraphics[width=0.3\textwidth]{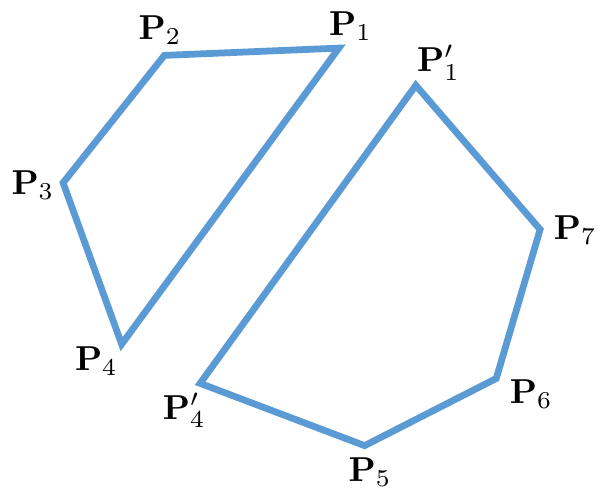} &
\includegraphics[width=0.3\textwidth]{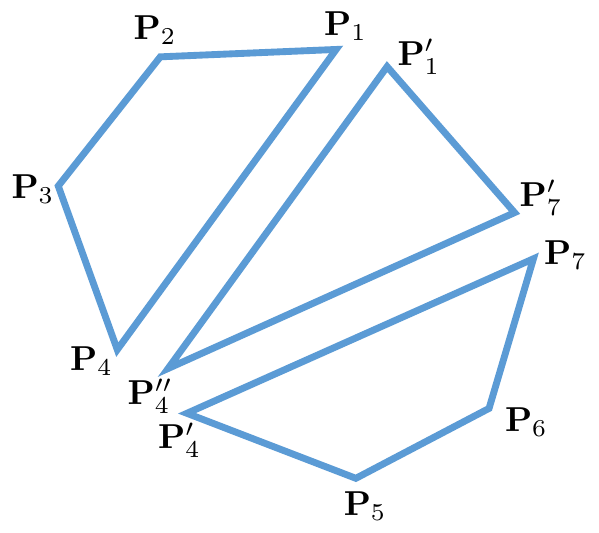} \\
(a) & (b) & (c)\\
\end{tabular}
\caption{Recursively splitting a polygon into quadrilaterals and triangles.}
\label{fig:poly_split}
\end{figure}

On the other hand, when a polygon has more than 4 edges (e.g., Fig.\ \ref{fig:poly_split}(a)), we split the polygon into two sub-polygons, one with 4 edges and the other with $n-2$ edges, where $n$ is the number of edges of the given polygon.
Two new edges are added due to splitting; see $\mathbf{P}_4 \mathbf{P}_1$ and $\mathbf{P}_1' \mathbf{P}_4'$ in Fig.\ \ref{fig:poly_split}(b).
Any splitting is valid as long as a quadrilateral is obtained.
We repeat this procedure for each resulting sub-polygon with more than 4 edges; see Fig.\ \ref{fig:poly_split}(c).
As a result, all resulting sub-polygons are either quadrilaterals or triangles, and a Coons patch is constructed accordingly for each of them as a cell.
This way, the polygonal domain is decomposed into the cells $\Tc_\ell$ ($\ell=1,\ldots,L$) whose associated mapping $\fundef{\Tb_\ell}{[0,1]^2}{\Tc_\ell}$ is a (degenerate) Coons patch.

\begin{figure}[htb]
\centering
\includegraphics[width=0.8\textwidth]{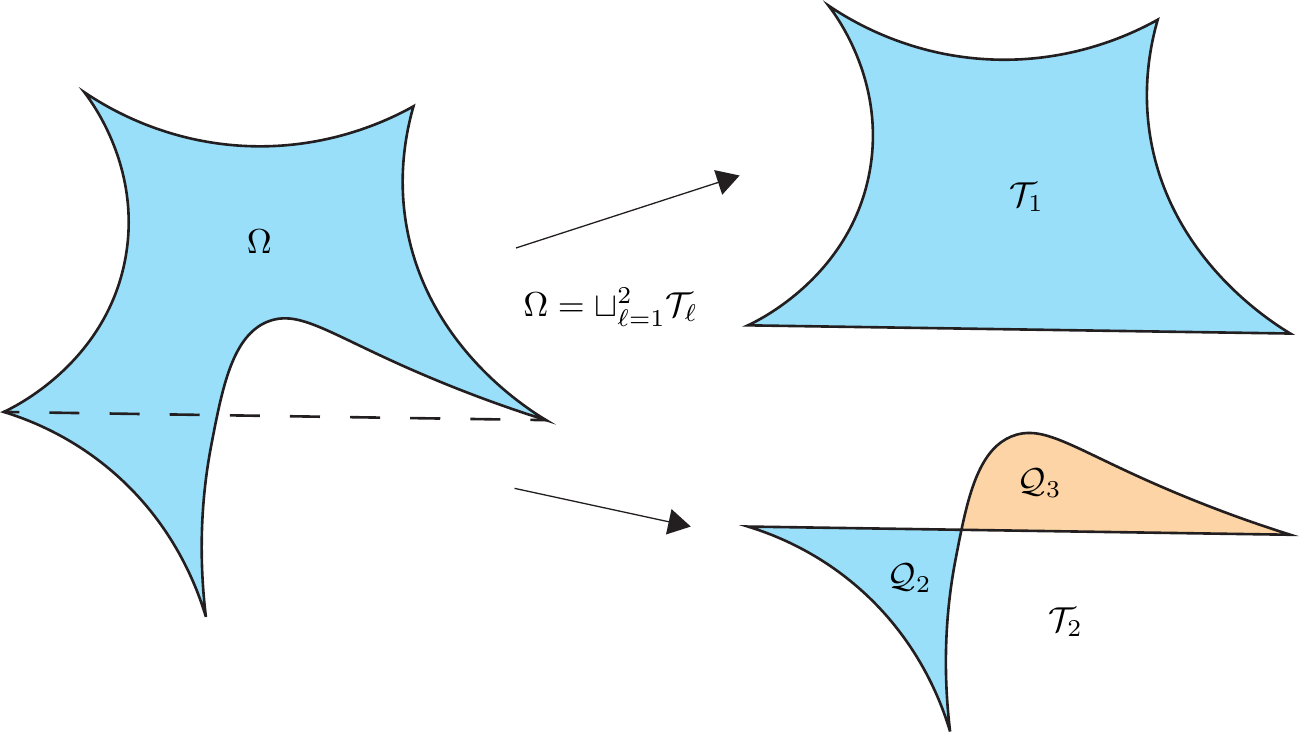}
\caption{Example of a folded decomposition such that $\Omega$ is split into $\{\Tc_1,\Tc_2\}$, where $\Tc_1$ is a $J^+$ cell and $\Tc_2$ is a folded cell.
$\Tc_2$ can be further split into $\{\Qc_2,\Qc_3\}$, where $\Qc_2$ is a $J^+$ cell and $\Qc_3$ is a $J^-$ cell.}
\label{fig:negative_poly}
\end{figure}

Edges of a curved polygonal domain generally correspond to high-order curves, and thus, an arbitrary decomposition may lead to folded cells.
For example, the decomposition in Fig.\ \ref{fig:negative_poly} is deemed valid, where, however, the splitting edge (dashed) intersects a polygonal edge and thus one of the resulting cells is a folded cell that has negative Jacobian in the orange region. We will show that such cells are valid cells for the integration purpose.

\subsection{Triangulation}
\label{sec:tri}

We create 2D triangulations to randomly generate folded decompositions to test the robustness of our proposed method.
In this paper, every triangle is treated as a Coons patch for simplicity.
In principle, any triangulation is allowed for a given polygon whose edges are counterclockwise oriented.
More specifically, a seed vertex $\mathbf{V}$ is first arbitrarily chosen in the bounding box of the polygon.
A triangle is then created with every polygonal edge and $\mathbf{V}$.

\begin{figure}[htb]
\centering
\begin{tabular}{cccc}
\includegraphics[width=0.2\textwidth]{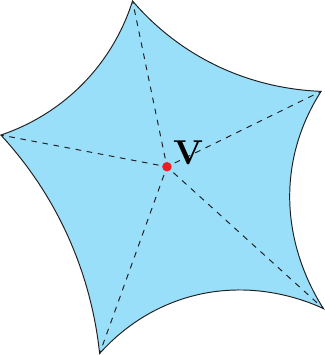} & 
\includegraphics[width=0.2\textwidth]{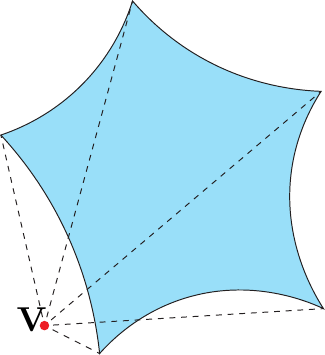} &
\includegraphics[width=0.3\textwidth]{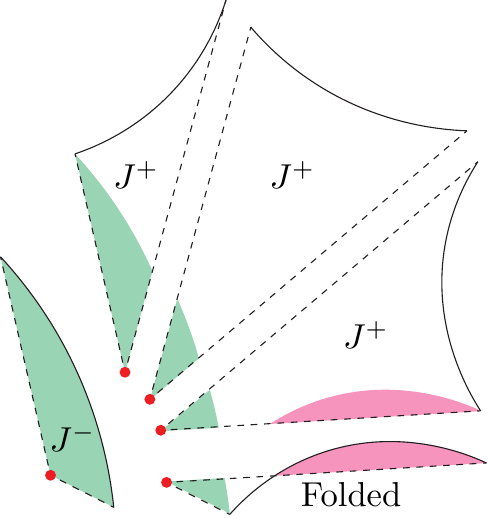} &
\includegraphics[width=0.2\textwidth]{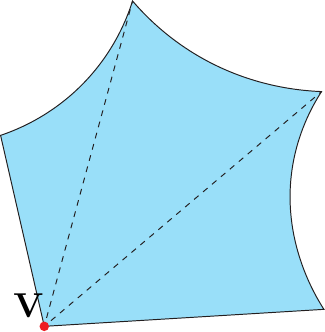} \\
(a) & (b) & (c) & (d)\\
\end{tabular}
\caption{Triangulations of curved polygons with different seed vertices $\mathbf{V}$.
(a) $\mathbf{V}$ inside the polygon, (b) $\mathbf{V}$ outside the polygon, (c) exploded view of (b), and (d) $\mathbf{V}$ being a vertex of the polygon.
In (c), the regions of the same color overlap with another.
The boundary curves of the polygons are oriented counterclockwise.}
\label{fig:tri}
\end{figure}

$\mathbf{V}$ can be inside or outside the polygon; see the red vertices in Figs.\ \ref{fig:tri}(a, b).
When $\mathbf{V}$ is inside, the generated decomposition may or may not have negative Jacobian, depending on whether or not a splitting edge intersects with any of the polygonal edges.
However, when $\mathbf{V}$ is outside, some of the resulting cells must have negative Jabobian; see Fig.\ \ref{fig:tri}(c).

Among all the choices of $\mathbf{V}$, a straightforward and sometimes efficient way is to pick a vertex of the polygon; see Fig.\ \ref{fig:tri}(d).
When an edge that shares this vertex is simply a straight segment, the resulting triangle has a zero area and thus will have no contribution to integration.
In other words, fewer cells are generated in this case; compare Fig.\ \ref{fig:tri}(c) with Figs.\ \ref{fig:tri}(a, b).

\begin{myremark}
While up to now we have assumed that the domain of interest is a planar domain, i.e., $\Omega\subset \R^2$, the proposed method also works for 2-manifold as long as it has a parametric representation.
More specifically, $\Omega$ can be a parameterized surface of the form, $\fundef{\mathbf{S}}{\mathcal{D}\subseteq [0,1]^2}{\R^3}$.
When this is the case, the proposed decomposition is applied to $\mathcal{D}$ rather than to $\Omega$.
As a result, $\mathcal{D}$ is decomposed into (possibly folded) cells $\Tc_\ell$ ($\ell=1,\ldots,L$), whose associated mappings $\fundef{\mathbf{T}_\ell}{[0,1]^2}{\Tc_\ell}$ are (degenerate) Coons patches.
\end{myremark}

\subsection{Theoretical study on folded decompositions}
\label{sec:theory}

In this section, we start with an intuitive explanation to motivate how folded decompositions work for integration, and then we provide a formal argument.

We first consider to compute the area of a curved polygon by integrating the Jacobian over all the cells that are obtained through either quadrilateral-dominant decomposition or triangulation.
We take decompositions in Figs.~(\ref{fig:poly_split}, \ref{fig:tri}) as examples.
In the case of Fig.\ \ref{fig:poly_split}(c) or Fig.\ \ref{fig:tri}(a) where all the cells are $J^+$ cells, it is clear that integrating the Jacobian over such cells leads to the area of the polygon.
On the other hand, when folded cells are involved, for example in Fig.\ \ref{fig:negative_poly} or Fig.\ \ref{fig:tri}(c), negative-Jacobian regions there overlap with certain positive-Jacobian regions in other cells, and thus their contributions are canceled out, leading to the desired area of the polygon; see the regions of matching colors in Fig.\ \ref{fig:tri}(c).
This implies that when the integrand is a constant, we can use folded decompositions to appropriately integrate it over a complex domain $\Omega$.

In fact, folded decompositions can be applied to much more general functions, which we will show in Section \ref{sec:int_poly} through numerical examples.
It is well known that canceling of positive and negative contributions may cause large deviations due to the round-off error.
On the other hand, we never spot such an issue in any of our examples.

\begin{myremark}
We note that for the cancellation to work, the boundary curves need to be oriented in a consistent manner (e.g., counterclockwise).
Moreover, we notice that due to the presence of negative Jacobian, generally a folded decomposition leads to a set of cells where $\Omega\subseteq \cup_{\ell=1}^L \Tc_\ell$. 
In other words, the resulting cells may cover a larger region than the domain of interest, which, therefore, requires the integrand $f(\xb)$ to be well defined on the entire $\cup_{\ell=1}^L \Tc_\ell$. 
We will discuss more on the definition of $f(\xb)$ in Section \ref{sec:iga}.
\end{myremark}

Next, we theoretically study folded decompositions for integration under the assumption that first, the integrand in Eq.~\eqref{eq:integral} is suitably defined over all the involved integration cells, and second, the antiderivative of the integrand exists.
The key idea is to apply the divergence theorem to transform an integral to its boundary counterpart.
We will show that the resulting boundary integral stays the same regardless of how the decomposition is achieved.
We work with the 2D case for illustration and explanation, but the argument holds for both 2D and 3D.
The following derivation is constructive rather than being comprehensive to provide insights on how folded decompositions work.

\begin{figure}[htb]
\centering
\begin{tabular}{ccc}
\includegraphics[width=0.3\textwidth]{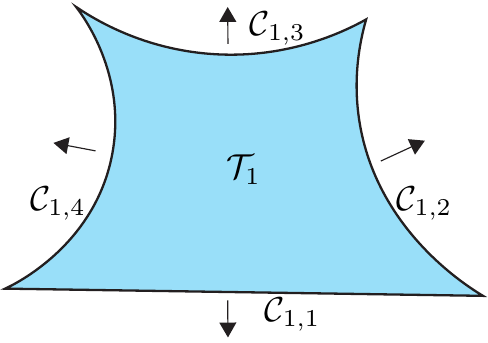} &\hspace{+2mm}
\includegraphics[width=0.35\textwidth]{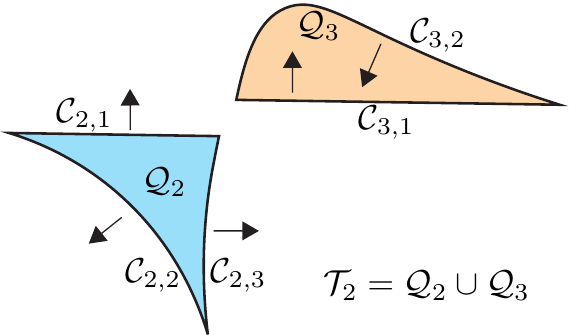} &\hspace{+2mm}
\includegraphics[width=0.22\textwidth]{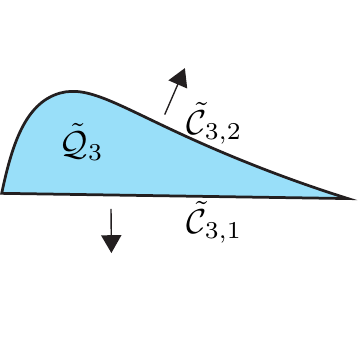} \\
(a) & (b) & (c)
\end{tabular}
\caption{The decomposed cells shown in Fig.\ \ref{fig:negative_poly} and notations of their boundary curves. (a) The $J^+$ cell, (b) the folded cell that consists of a $J^+$ subcell and a $J^-$ subcell, and (c) reparameterization of $\Omega_3$ to obtain a $J^+$ counterpart.
The black arrows indicate the normal directions.
}
\label{fig:negative_div}
\end{figure}

In particular, we consider the domain $\Omega$ shown in Fig.\ \ref{fig:negative_poly}, where the decomposition leads to a $J^+$ cell $\Tc_1$ and a folded cell $\Tc_2$. 
Initially, the integral of interest is computed using $\Tc_1$ and $\Tc_2$, 
\begin{equation}
\int_{\Omega} f(\xb) = \int_{\Tc_1} f(\xb) +  \int_{\Tc_2} f(\xb).
\label{eq:div_1}
\end{equation}
We further split $\Tc_2$ into $\Tc_2=\Qc_2\cup\Qc_3$ such that $\Qc_2$ is a $J^+$ cell and $\Qc_3$ is a $J^-$ cell; see Fig.~\ref{fig:negative_div}.
Specifically, letting $\fundef{\Qb_k}{[0,1]^2}{\Qc_k}$ ($k=2,3$), we have
\begin{equation}
\int_{\Tc_2} f(\xb) = \sum_{k=2}^3 \int_{\Qc_k} f(\xb) = \sum_{k=2}^3  \int_{[0,1]^2} f(\Qb_k(\ub)) \jacob{\Qb_k(\ub)},
\end{equation}
where $\jacob{\Qb_3(\ub)} \leq 0$.
To attain cells of consistent Jacobian, we reparameterize $\Qc_3$ as $\tilde{\Qc}_3$ such that 
\begin{equation}
\int_{\tilde{\Qc}_3} f(\xb) = \int_{[0,1]^2} f(\tilde{\Qb}_3(\mathbf{s})) \jacob{\tilde{\Qb}_3(\mathbf{s})} =  -\int_{[0,1]^2} f(\Qb_3(\ub)) \jacob{\Qb_3(\ub)} = - \int_{\Qc_3} f(\xb),
\end{equation}
where $\tilde{\Qc}_3=\Qc_3$, $\tilde{\Qb}_3$ is the mapping of $\tilde{\Qc}_3$, and $\jacob{\tilde{\Qb}_3} \geq 0$.
Eq. \eqref{eq:div_1} becomes
\begin{equation}
\int_{\Omega} f(\mathbf{x}) = \int_{\Tc_1} f(\mathbf{x}) +  \int_{\Qc_2} f(\mathbf{x}) -  \int_{\tilde{\Qc}_3} f(\mathbf{x}).
\label{eq:div_2}
\end{equation}

Now we assume the antiderivative of $f(\mathbf{x})$ to be $\mathbf{F}(\mathbf{x})$ such that $f(\mathbf{x})=\nabla \cdot \mathbf{F}(\mathbf{x})$.
Applying the divergence theorem, we have
\begin{equation}
\begin{aligned}
\int_{\Omega} f(\mathbf{x}) &= \int_{\partial\Tc_1} \mathbf{F}(\mathbf{x}) \cdot \mathbf{n}_1 +  \int_{\partial\Qc_2} \mathbf{F}(\mathbf{x}) \cdot \mathbf{n}_2 -  \int_{\partial\tilde{\Qc}_3} \mathbf{F}(\mathbf{x}) \cdot \tilde{\mathbf{n}}_3 ,\\
&= \sum_{j=1}^4 \int_{\mathcal{C}_{1,j}} \mathbf{F}(\mathbf{x}) \cdot \mathbf{n}_1 +  \sum_{j=1}^3 \int_{\mathcal{C}_{2,j}} \mathbf{F}(\mathbf{x}) \cdot \mathbf{n}_2 - \sum_{j=1}^2 \int_{\tilde{\mathcal{C}}_{3,j}} \mathbf{F}(\mathbf{x}) \cdot \tilde{\mathbf{n}}_3,\\
\end{aligned}
\label{eq:div_2bis}
\end{equation}
where $\mathbf{n}_1$, $\mathbf{n}_2$, and $\tilde{\mathbf{n}}_3$ are the outwards normals to corresponding domains, and $\mathcal{C}_{i,j}$ are the boundary curves; see Fig.~\ref{fig:negative_div}.
Note that we have $\mathbf{n}_2 = -\mathbf{n}_1$ when $\mathbf{x}\in\mathcal{C}_{2,1}$, and $\tilde{\mathbf{n}}_3 = \mathbf{n}_1$ when $\mathbf{x}\in\tilde{\mathcal{C}}_{3,1}$.
Therefore, contributions from common curves (i.e., $\mathcal{C}_{2,1}$ and $\tilde{\mathcal{C}}_{3,1}$ with $\mathcal{C}_{1,1}$) are canceled out, leading to
\begin{equation}
\int_{\Omega} f(\mathbf{x}) = \sum_{j=2}^4\int_{\mathcal{C}_{1,j}} \mathbf{F}(\mathbf{x}) \cdot \mathbf{n}_1 + \sum_{j=2}^3 \int_{\mathcal{C}_{2,j}} \mathbf{F}(\mathbf{x}) \cdot \mathbf{n}_2 - \int_{\tilde{\mathcal{C}}_{3,2}} \mathbf{F}(\mathbf{x}) \cdot \tilde{\mathbf{n}}_3.
\label{eq:div_3}
\end{equation}
Regarding the last term in Eq.~\eqref{eq:div_3}, after reversing the direction of $\tilde{\mathcal{C}}_{3,2}$, we have
\begin{equation}
\int_{\tilde{\mathcal{C}}_{3,2}} \mathbf{F}(\mathbf{x}) \cdot \tilde{\mathbf{n}}_3 = -\int_{\mathcal{C}_{3,2}} \mathbf{F}(\mathbf{x}) \cdot \mathbf{n}_3,
\label{eq:div_4}
\end{equation}
where $\mathbf{n}_3 = -\tilde{\mathbf{n}}_3$.
$\mathcal{C}_{3,2}$ has the same orientation as $\partial\Omega$ and $\mathbf{n}_3$ is the outwards normal to $\partial\Omega$ when $\mathbf{x}\in\mathcal{C}_{3,2}$.
Therefore, examining Eqs.~(\ref{eq:div_3}, \ref{eq:div_4}) reveals that the integral over $\Omega$ reduces to the boundary integral,
\begin{equation}
\int_{\Omega} f(\mathbf{x}) = \int_{\partial\Omega} \mathbf{F}(\mathbf{x}) \cdot \mathbf{n},
\label{eq:div_5}
\end{equation}
where we note that $\partial\Omega=\mathcal{C}_{1,2}\cup\mathcal{C}_{1,3}\cup\mathcal{C}_{1,4}\cup\mathcal{C}_{2,2}\cup\mathcal{C}_{2,3}\cup\mathcal{C}_{3,2}$.

Therefore, Eqs.~(\ref{eq:div_1}, \ref{eq:div_5}) are equivalent.
As the involved decomposition represents a generic case and the same argument can be straightforwardly applied to other decompositions, we conclude that the integral of interest leads to the same boundary integral regardless of how a decomposition is created, given that the above-mentioned two assumptions hold.
In other words, if a $J^+$ decomposition exists for the domain shown in Fig.\ \ref{fig:negative_poly}, then theoretically it leads to the same integration result as that using the folded decomposition.
Indeed, we will show different kinds of numerical evidence later.

\section{Folded decompositions for $d=3$}
\label{sec:nj3d}

Now we discuss the 3D case, that is, to use folded decompositions for integration over curved polyhedra.
The idea is analogous to Section \ref{sec:tri}, but here our goal is to decompose a given polyhedron into a set of (triangular and/or rectangular) pyramids, which, in particular, can naturally deal with gaps/overlaps as well as sharp features in the input B-rep.
Pyramids are treated as degenerated hexahedra.
The decomposition needs two steps: (1) decomposition of the input B-rep following the previous section, and (2) search for a seed vertex to build pyramids.

As the input B-rep is essentially a 2-manifold, Step (1) immediately follows Section~\ref{sec:nj2d}. We assume that the B-rep has the form $\po=\{\mathbf{S}_i(u,v)\in\R^3,\, (u,v)\in\mathcal{D}_i\subseteq [0,1]^2,\,i=1,\ldots,N_f\}$. 
It is worth mentioning that the mesh is not required to be conformal on the common edge of any two faces. 
Such edges sometimes represent sharp features of $\po$. 
Moreover, gaps/overlaps are allowed in the input B-rep as long as the topological information is correct. 
Such settings are supported because folded decompositions are generated independently for each $\mathbf{S}_i(u,v)$, or more precisely, for $\mathcal{D}_i$ in the parametric domain, leading to $\fundef{\mathbf{R}_{i,j}}{[0,1]^2}{\R^2}$, $j=1,\ldots,N_i$, where $\mathbf{R}_{i,j}$ is a (degenerate) Coons patch that may have negative Jacobian.
As a result, each face $\mathbf{S}_i$ is decomposed as $\mathbf{S}_i\circ\mathbf{R}_{i,j}([0,1]^2)$, $j=1,\ldots,N_i$.
Each $\mathbf{S}_i\circ\mathbf{R}_{i,j}([0,1]^2)$ is further approximated by a B\'{e}zier surface to obtain a polynomial representation.



Next in Step (2), we choose a seed vertex $\mathbf{V}$ to create target pyramids with all the resulting patches $\mathbf{S}_i\circ\mathbf{R}_{i,j}([0,1]^2)$ from Step (1), where $\mathbf{V}$ and $\mathbf{S}_i\circ\mathbf{R}_{i,j}$ serve as the top vertex and the base surface, respectively.
While $\mathbf{V}$ can be chosen in a rather arbitrary manner, we usually choose it to be one of the vertices of the polyhedron as it can yield fewer cells (similar to 2D).

Once $\mathbf{V}$ and $\mathbf{R}_{i,j}(u,v)$ are in place, a pyramidal cell $\mathbf{T}_{i,j}(u,v,w)$ is created as a ruled volume between $\mathbf{V}$ and each $\mathbf{R}_{i,j}(u,v)$,
\begin{equation}
\mathbf{T}_{i,j}(u,v,w) = (1-w) \mathbf{S}_i \circ \mathbf{R}_{i,j}(u,v) + w \mathbf{V}, \quad (u,v,w)\in [0,1]^3.
\label{eq:cell3d}
\end{equation}
Therefore, the given curved polyhedron is decomposed into a set of such cells $\{\mathbf{T}_{i,j}(u,v,w)\}$ for integration.

\begin{figure}[htb]
\centering
\begin{tabular}{ccc}
\includegraphics[width=0.3\textwidth]{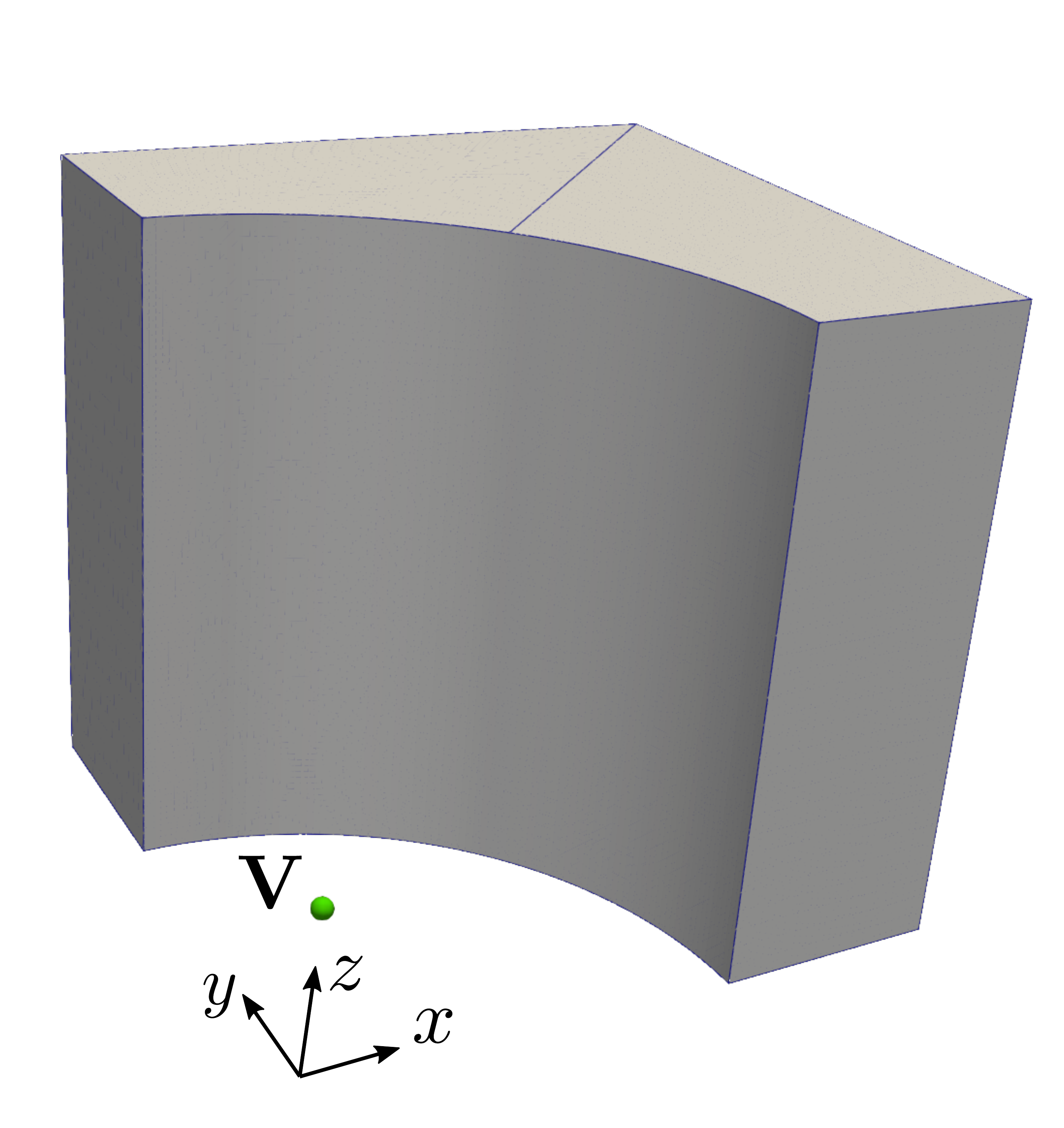} & 
\includegraphics[width=0.3\textwidth]{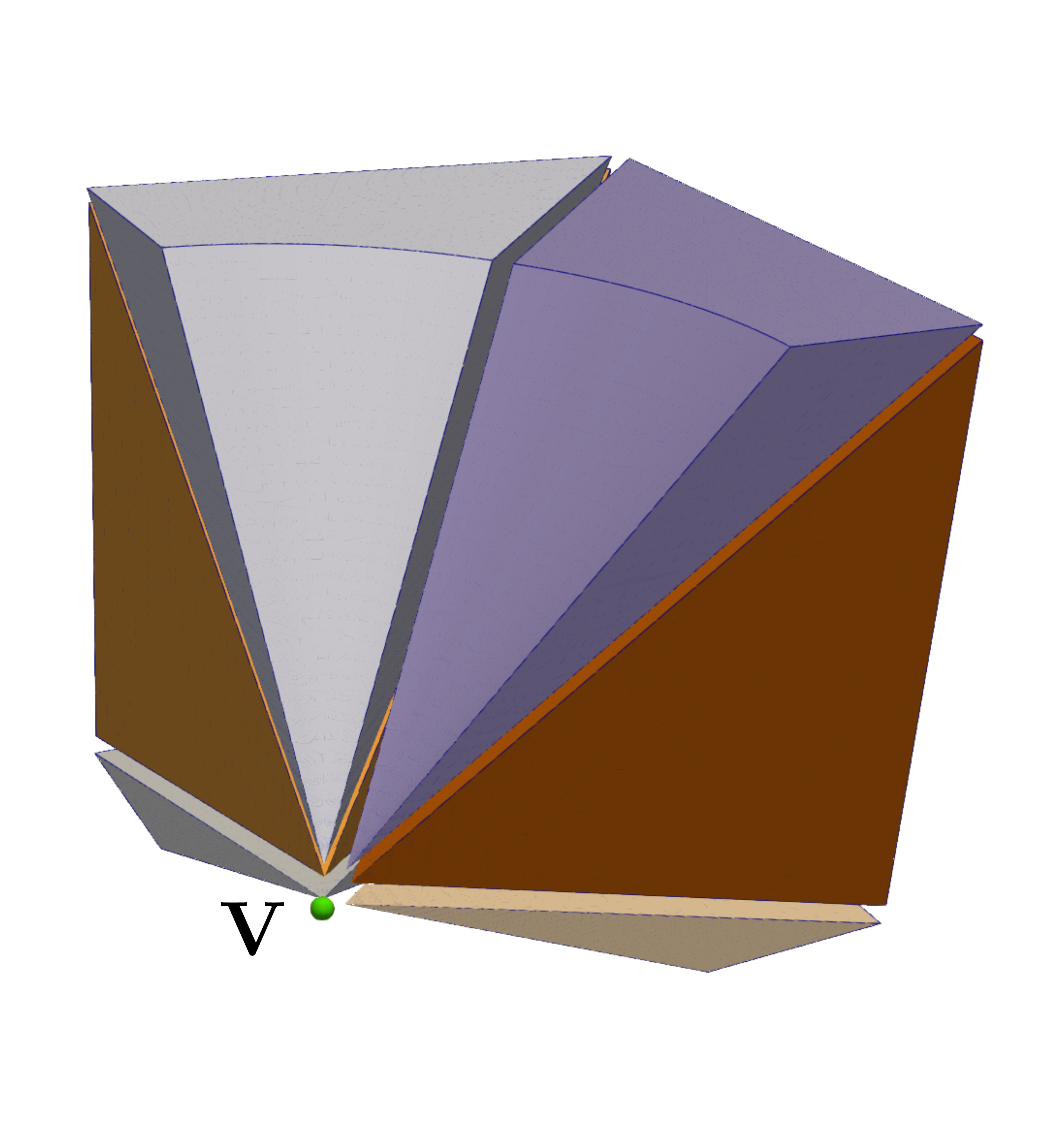} & 
\includegraphics[width=0.3\textwidth]{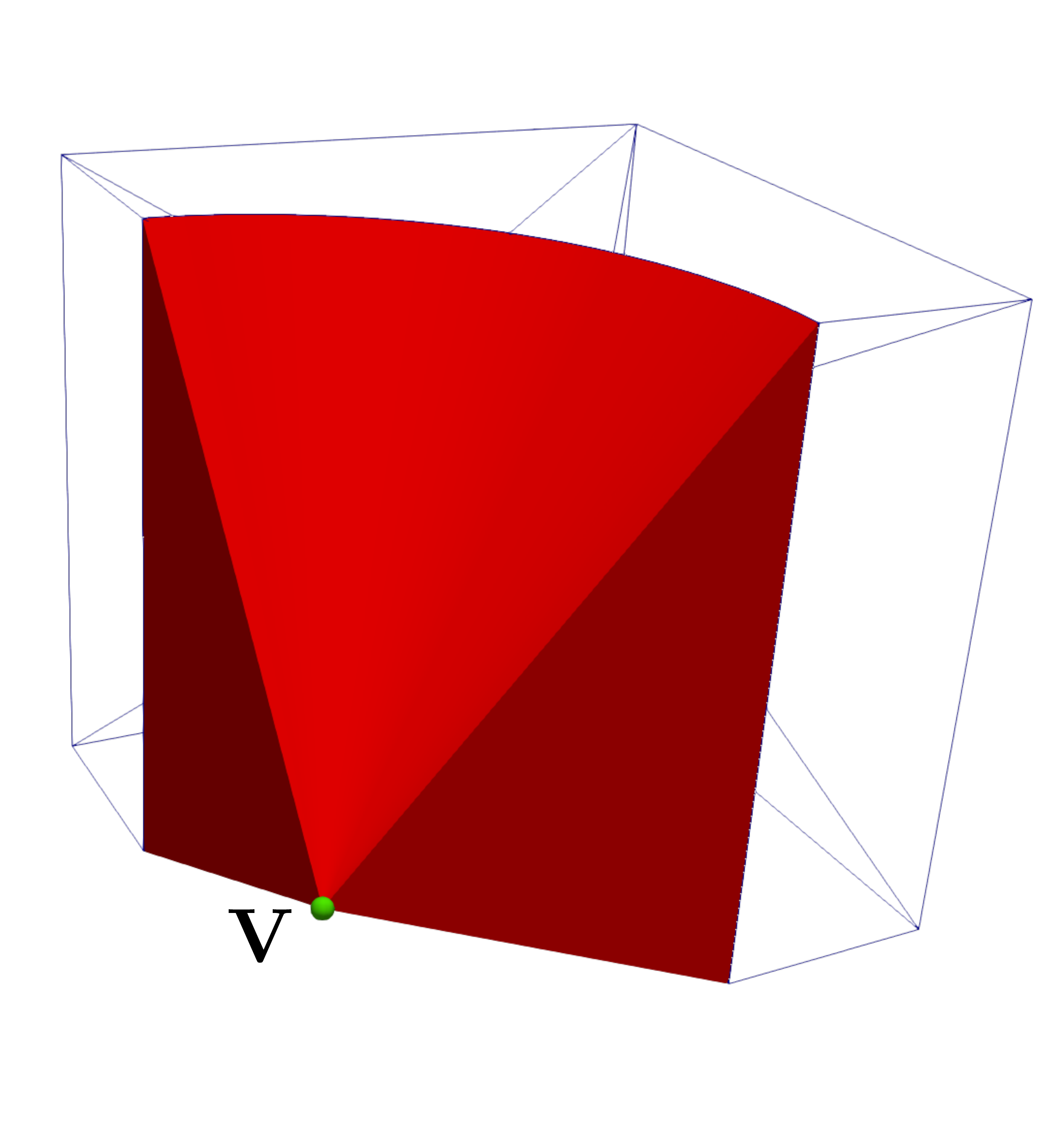} \\
(a) & (b) & (c)\\
\end{tabular}
\caption{Split of a curved polyhedron using an arbitrary vertex $\mathbf{V}$.
(a) The input B-rep (a cube cut with a cylinder); (b) six of the resulting cells present positive Jacobian; and (c) one of the cells has negative Jacobian.}
\label{fig:poly_decomp_3D}
\end{figure}

As an example, in Fig.\ \ref{fig:poly_decomp_3D} the decomposition of a 3D B-rep is shown.
The polyhedra to decompose is a cube $[0,1]^3$ cut with a cylinder of radius $0.65$ and axis along the $z$ direction (see Fig.~\ref{fig:poly_decomp_3D}(a)).
A seed vertex $\mathbf{V}$ is (arbitrarily) chosen to be at $(0,0,0.3)$.
The bottom and top faces of the B-rep (5-sided polygons) are decomposed into two quadrilaterals each, whereas all the other faces require a single quadrilateral.
By generating pyramids between the generated surface patches and the vertex $\mathbf{V}$ a total of 7 cells are created.
Six of them (Fig.\ \ref{fig:poly_decomp_3D}(b)) present positive Jacobian everywhere except at the vertex $\mathbf{V}$, where the Jacobian vanishes, so they are $J^+$ cells.
On the other hand, the cell shown in Fig.\ \ref{fig:poly_decomp_3D}(c) has negative Jacobian everywhere except at $\mathbf{V}$, where the Jacobian vanishes, so it is a $J^-$ cell.

\begin{myremark}
The work \cite{ref:massarwi19} strives to find a so-called kernel point as the seed vertex such that all the resulting cells can have positive Jacobian.
A kernel point is a point that is visible to all the points on the boundary.
This is known as the three-dimensional \emph{art gallery} problem.
Generally, it involves recursively subdividing the given domain to find an appropriate kernel point for each subdomain, which complicates implementation and generates a large number of cells.
Due to the great flexibility folded decompositions provide, we do not need a kernel point for decomposition in this work.
\end{myremark}

\section{Application: immersed methods}
\label{sec:iga}

One natural application of the integration method described above is isogeometric analysis (IGA) \cite{ref:hughes05} with trimmed volumetric geometries.
The proposed method can be also used in other applications such as virtual element methods \cite{ref:veiga13}.
One typical way of using such representations in IGA is through immersed methods \cite{ref:schillinger12}, which embed a B-rep into a background Cartesian grid and thus leads to a non-conformal description of the computational domain; see Fig.\ \ref{fig:ex_poly_curve}(a) for example.

Numerical integration in IGA is performed at the element level, so we first classify different types of elements.
Let $\mathcal{G}$ be a Cartesian grid, and $K\in\mathcal{G}$ be an element.
The computational domain is $\Omega$ and $\Gamma=\po$.
Depending on how $K$ relates to $\Gamma$ and $\Omega$, elements can be divided into trimmed elements ($K\cap \Gamma\neq \varnothing$), non-trimmed active elements ($K\cap \Omega= K$), and non-trimmed inactive elements ($K\cap \Omega= \varnothing$) that are not assembled in general.
While it is straightforward to use standard Gauss quadrature for non-trimmed elements, trimmed elements need special treatment.

Every trimmed element ($K_a=K\cap\Omega\neq\varnothing$) is a curved polygon ($d=2$) or a curved polyhedron ($d=3$), to which we apply our decomposition strategy.
In other words, $K_a$ is now the domain of interest in Eq.~\eqref{eq:integral}.
The boundary of $K_a$, e.g., $K\cap\Gamma$ in particular, generally does not have a polynomial representation, and instead it has a piecewise polynomial structure.
To build integration cells, polynomial fitting is needed to find a corresponding \bz representation to approximate $K\cap\Gamma$.
The degree of the \bz representation is left as a parameter to be chosen to guarantee the desired level of accuracy; see, e.g.,~\cite{ref:antolin19a}.




\begin{myremark}
Recall that the tensor-product Gauss quadrature is defined on the parent domain $\Lambda=[0,1]^d$ ($d=2,3$) of each cell.
The positions of quadrature points in a trimmed element $K_a$ (whose parent element is $K$) are obtained by the geometric mapping $\fundef{\Tb_{K_a}}{[0,1]^d}{\Tc}$ of the cell; see Eqs.\ (\ref{eq:cell2d}, \ref{eq:cell3d}).
Their quadrature weights in the domain of $K_a$ are the product of the Jacobian of $\mathbf{T}$ with the corresponding weights of Gauss quadrature in $\Lambda$; see also Eq.\ \eqref{eq:gptw}.
This implies that the weights can be negative or zero, which, however, is not an issue because they are only used in multiplication and never in division.
However, due to the nature of negative Jacobian, we may have $\mathbf{T}(\Lambda)\not\subset K_a$, and thus some quadrature points may lie outside of $K_a$. 
When certain quantities of interest are not well defined outside $K_a$, one solution to this potential issue is through polynomial extension \cite{ref:puppi20}, which enables a stable way to evaluate quantities outside of the domain.
\end{myremark}

\begin{myremark}
Recall that when creating integration cells, we do not require their parameterization to be conformal to their common interfaces.
This not only provides great flexibility on how we can decompose a given domain, but also supports B-reps that have artificial defects such as gaps and overlaps, thus opening the door to deal with dirty CAD geometries in a robust way.
\end{myremark}

\begin{myremark}
While folded decompositions provide a significant flexibility in accommodating real-world applications, there are limitations due to their nature of exhibiting negative Jacobian.
For instance, visualizing simulation results with folded cells poses an issue.
The sign of Jacobian can be seen as the orientation of a local region.
Change of sign means change of orientation, and thus visualization detects a reverse reflection on negative-Jacobian regions.
As a result, negative-Jacobian regions are invisible to users.
The compromise we make in this work is to create $J^+$ cells for visualization wherever it is needed.
For example, the boundary of a volume is decomposed according to quadrangulation \cite{ref:wei21b} to ensure that all the surface cells have positive Jacobians, where quantities of interest are computed and visualized.
This is done independently from creating volumetric cells, and it is only meant for visualization purposes.
\end{myremark}

\section{Integration of regular functions}
\label{sec:int_poly}

In this section, we demonstrate the efficacy and robustness of folded decompositions by integrating several regular functions.
We first compare the performance of folded and $J^+$ decompositions in computing integrals of polynomials.
Next, we study the convergence behaviors in integrating general functions that are not polynomials.

\subsection{Integration of polynomials}

The domain of interest is denoted as $\Omega$.
The boundary $\Gamma$ of $\Omega$ is embedded into a Cartesian grid $\mathcal{G}$, where we are only interested in trimmed elements as polynomials can be easily integrated exactly in non-trimmed elements.
For each trimmed element $K_a=K\cap \Omega\neq\varnothing$ ($K\in\mathcal{G}$ and $K\cap\Gamma\neq\varnothing$), we aim to compute
$$
I_{K,i} = \int_{K_a} B_i^p(\mathbf{x}),
$$
where $B_i^p(\mathbf{x})$ is $i$-th degree-$p$ tensor-product Bernstein polynomial defined in $K$.
Note that $\mathrm{span}\{B_i^p(\mathbf{x})\}_i = \Q_{p,p}(K)$ in 2D and $\mathrm{span}\{B_i^p(\mathbf{x})\}_i = \Q_{p,p,p}(K)$ in 3D.
Generally, $K_a$ is approximated with a set of integration cells $\{\mathbf{T}_{K,j}(\Lambda)\}_j$ ($\Lambda=[0,1]^d$), and thus $I_{K,i}$ is computed as
\begin{equation}
I_{K,i} = \int_{K_a} B_i^p(\mathbf{x}) \approx \sum_{j} \int_{\Lambda} B_i^p \circ \mathbf{T}_{K,j} \det \left( \nabla \mathbf{T}_{K,j} \right).
\label{eq:IKi}
\end{equation}
We use the tensor-product Gauss quadrature rule with $n$ points per direction in $\Lambda$.
An overkill result with 64 quadrature points per direction is taken as the reference ($I_{K,i}^{\mathrm{ref}}$) to evaluate error,
\begin{equation}
\mathrm{Err} = \max_{K\in\mathcal{G}, K\cap\Gamma \neq\varnothing} \max_{i} | I_{K,i}^h - I_{K,i}^{\mathrm{ref}} |,
\label{eq:err}
\end{equation}
where $I_{K,i}^h$ is the approximate to $I_{K,i}$ obtained by \eqref{eq:IKi} with Gauss quadrature.

\subsubsection{B-rep in 2D with a B-spline curve}
\label{sec:square_polynomial}

\begin{figure}[htb]
\centering
\begin{tabular}{cc}
\includegraphics[width=0.4\textwidth]{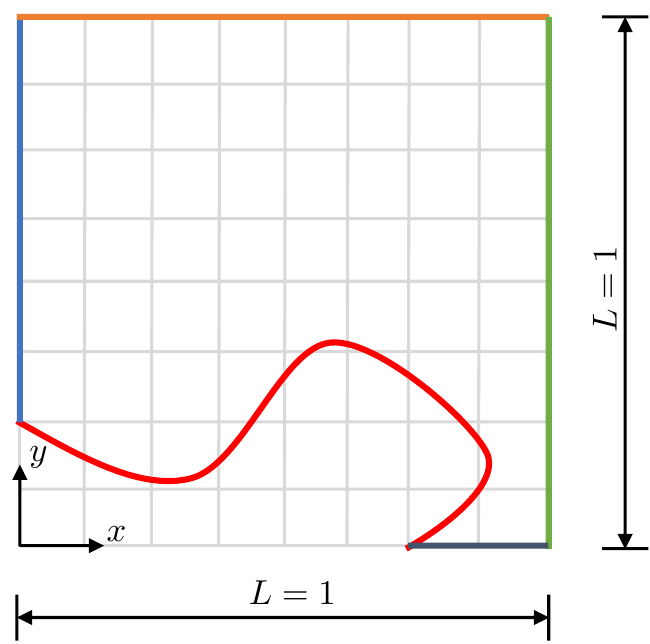} &\hspace{+2mm}
\includegraphics[width=0.47\textwidth]{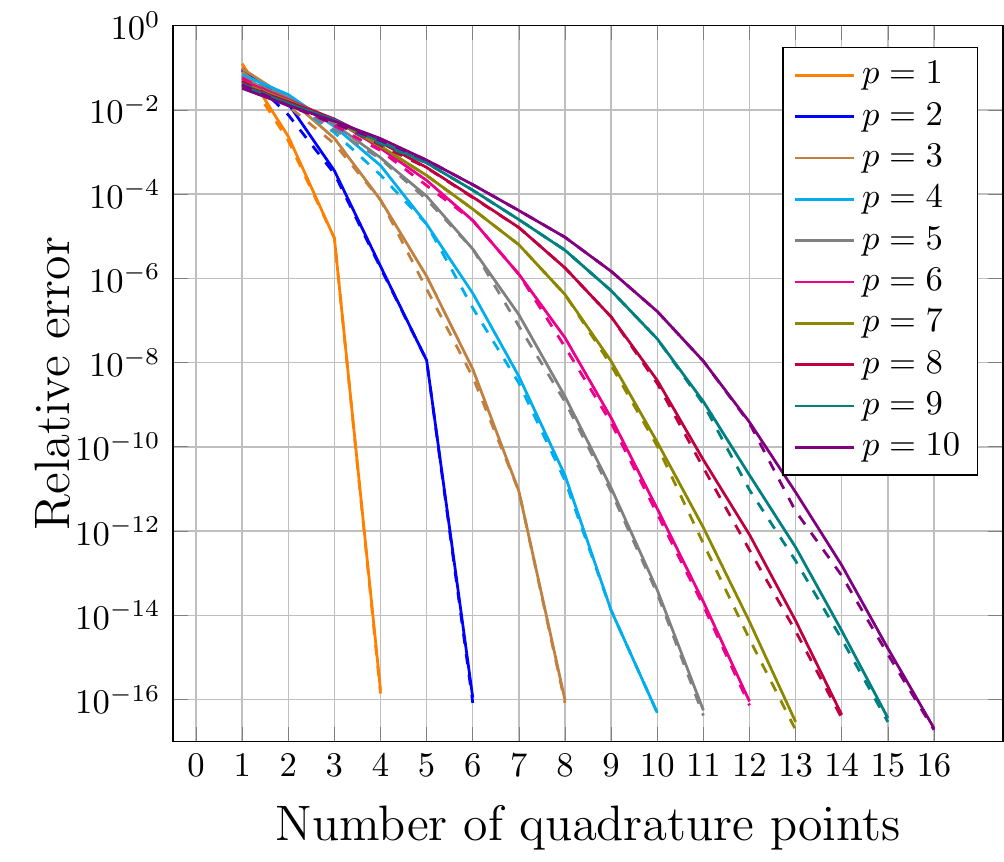} \\
(a) & (b)
\end{tabular}
\caption{Integration of polynomials on a 2D B-rep that involves a quadratic B-spline curve (red).
(a) The input B-rep composed of 5 connected curves, and (b) relative integration error with respect to the number of quadrature points per direction, where different colors represent different degrees, and solid and dashed curves are the results using folded and $J^+$ decompositions, respectively.}
\label{fig:ex_poly_curve}
\end{figure}

The geometry of the first 2D example is shown in Fig.\ \ref{fig:ex_poly_curve}(a).
Its B-rep consists of 5 connected curves: 4 straight segments and a quadratic B-spline curve (red), whose knot vector is $\{0,0,0,0.25,0.5,0.75,1,1,1\}$ and its control points are $(0,0.25)$, $(0.25,0)$, $(0.5,0.5)$, $(0.9,0.25)$, $(0.8,0.125)$, and $(0.75,0)$.
The B-rep is embedded into a $8\times 8$ Cartesian grid.
We consider two ways of decomposition for trimmed elements: the $J^+$ decomposition by quadrangulation~\cite{ref:wei21b}, and the folded decomposition through triangulation in Section \ref{sec:tri}.
In either case, the degree of all the \bz patches (i.e., $\mathbf{T}_{K,j}$ in \eqref{eq:IKi}) is set to be 2, the same as that of the B-spline curve.

In the case of folded decomposition, we choose different seed vertices to create a family of triangulations for each trimmed element $K_a$ (whose parent element is $K$).
Specifically, $11\times 11$ points are sampled according to a uniform Cartesian distribution in $K$, and each of them serves as a seed vertex to create a corresponding triangulation.
This is aimed to test whether folded decompositions can work robustly in a rather arbitrary setting.
Different triangulations may yield different integration errors, among which the largest is taken as the error of folded decomposition.

We use different degrees of Bernstein polynomials (ranging from 1 to 10) to test the integration accuracy with respect to the number of quadrature points per direction.
The results are summarized in Fig.\ \ref{fig:ex_poly_curve}(b).
We observe that the results are almost indistinguishable between folded and $J^+$ cells.
Moreover, according to \eqref{eq:n2d} in~\ref{appA}, where we study the number of quadrature points required to exactly integrate a polynomial, we need at least $2(p+1)$ quadrature points (per direction) to exactly integrate $B_i^p(\mathbf{x})$.
In particular, we observe this agreement in Fig.~\ref{fig:ex_poly_curve}(b) for degrees from 1 to 4, where the number of quadrature points needed to achieve machine precision ($\sim 10^{-16}$) is 4, 6, 8, and 10, respectively.
However, for even higher degrees, machine precision is reached before $2(p+1)$ quadrature points are used, and further increasing the number of quadrature points will no longer improve the accuracy.

\subsubsection{B-rep in 2D with a rational curve}

\begin{figure}[!htb]
\centering
\begin{tabular}{cc}
\includegraphics[width=0.4\textwidth]{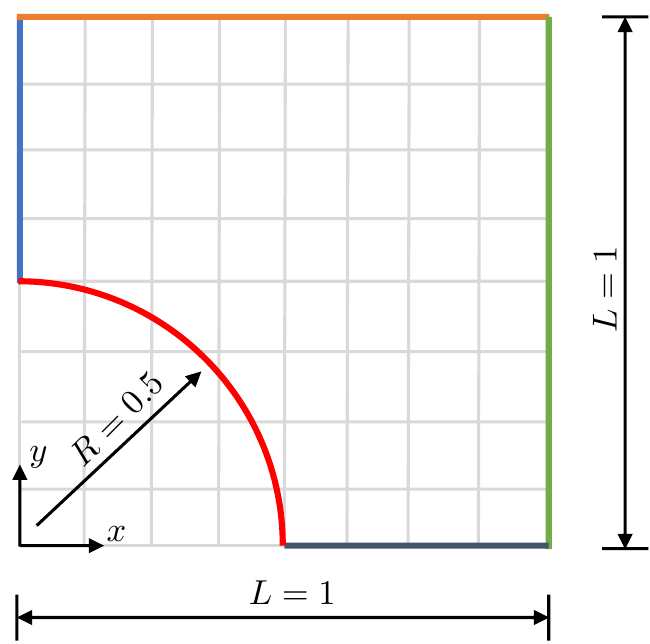} &\hspace{+2mm}
\includegraphics[width=0.47\textwidth]{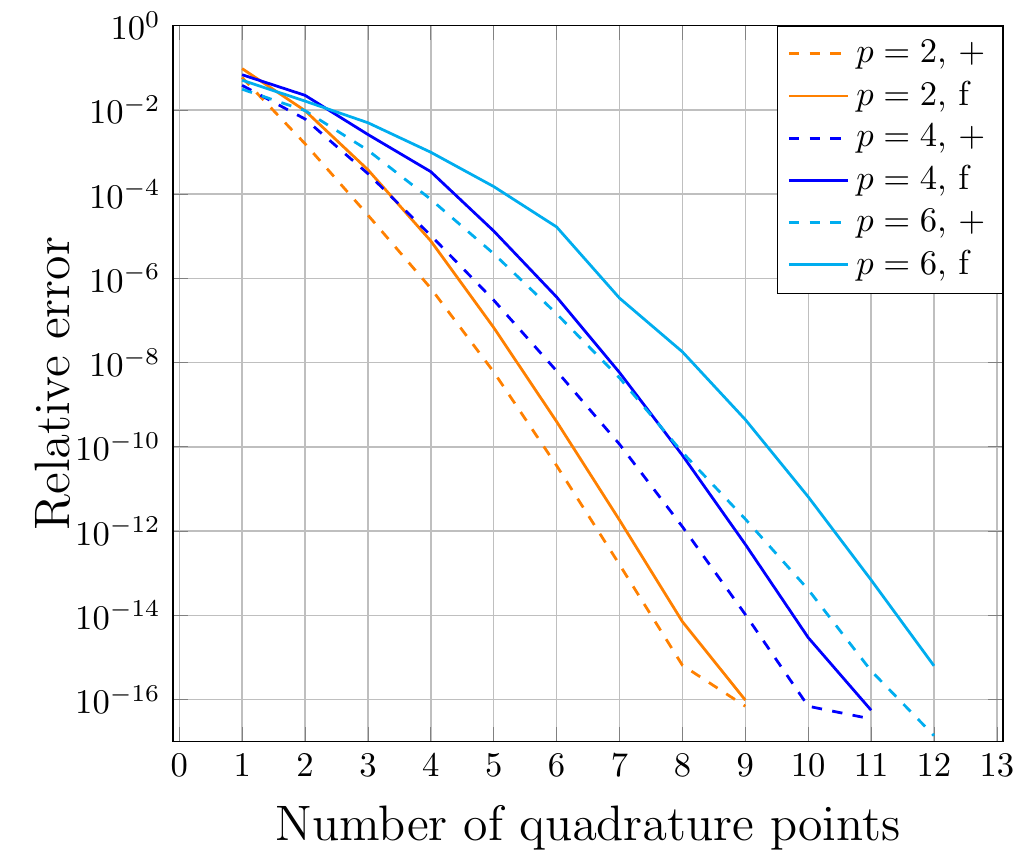} \\
(a) & (b)
\end{tabular}
\caption{Integration of polynomials on a 2D B-rep that involves a quadratic rational \bz curve (red).
(a) The input B-rep composed of 5 connected curves, and (b) relative integration error with respect to the number of quadrature points per direction, where different colors represent different degrees, and solid and dashed curves are the results using folded and $J^+$ decompositions, respectively.}
\label{fig:ex_ration}
\end{figure}

Next, we study another 2D B-rep that involves a rational curve.
The input is shown in Fig.\ \ref{fig:ex_ration}(a), where the red curve is a quadratic rational \bz curve exactly representing a quarter of a circle.
Other than the difference in the B-rep, everything else follows the same as that in Section \ref{sec:square_polynomial}.
The degrees of Bernstein polynomials for integration are chosen to be $p\in\{2,4,6\}$.
The circle is approximated with polynomial \bz curves of degree $p$, which are then used to create in either folded or $J^+$ decompositions.
The results are summarized in Fig.\ \ref{fig:ex_ration}(b).
We observe the following.
First, overall, $J^+$ decompositions yield more accurate results than folded decompositions, where fewer quadrature points are used to achieve the same accuracy.
However, recall that we create $11\times 11$ triangulations for each cut element, and we take the worst integration results to be the results using folded decompositions in Fig.\ \ref{fig:ex_ration}(b).
Second, with a fixed number of quadrature points, the relative accuracy difference between folded and $J^+$ decompositions is really not significant; compare solid and dashed curves when the relative error is $10^{-16}$.
Third, folded decompositions can also achieve machine precision by using slightly more quadrature points than $J^+$ decompositions.
Fourth, $J^+$ and folded decompositions show the same convergence rates, while the latter may have a longer pre-asymptotic regime.

\subsubsection{B-rep in 3D with a polynomial surface}
\label{sec:cube_polynomial}

\begin{figure}[htb]
\centering
\begin{tabular}{cc}
\includegraphics[width=0.47\textwidth]{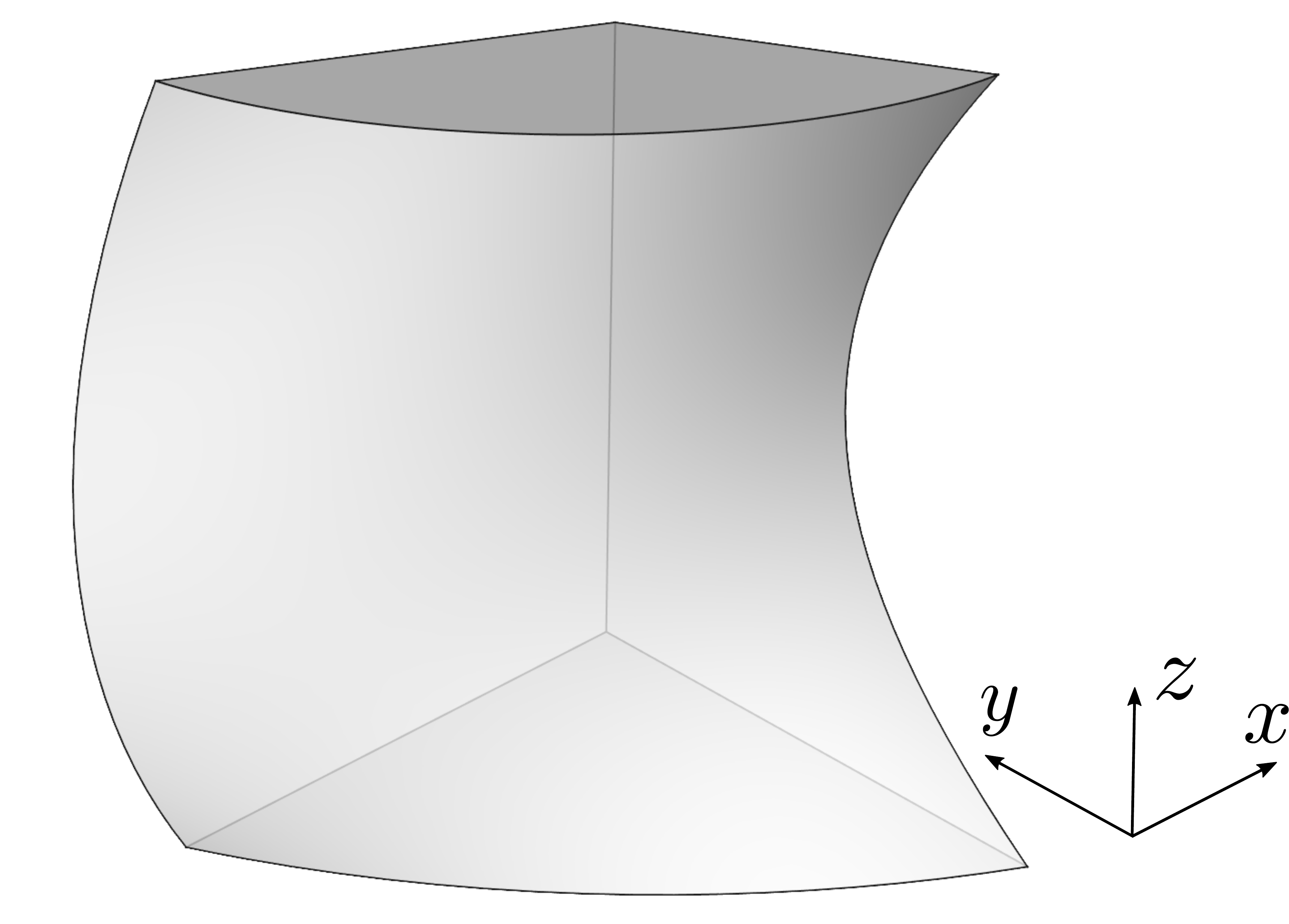} &\hspace{+2mm}
\includegraphics[width=0.47\textwidth]{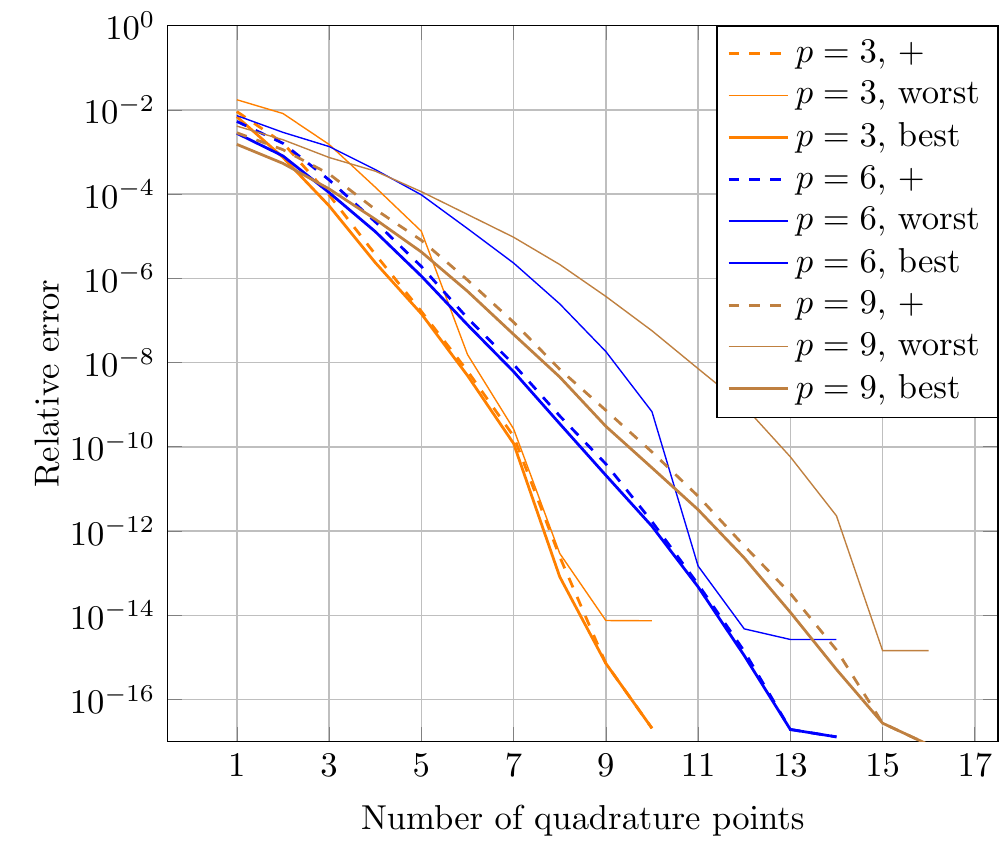} \\
(a) & (b)
\end{tabular}
\caption{Integration of polynomials on a 3D B-rep that involves a bi-quadratic \bz surface.
(a) The input B-rep composed of 4 axis-aligned planar surfaces and a bi-quadratic \bz surface, and (b) relative integration error with respect to the number of quadrature points per direction, where different colors represent different degrees, and solid and dashed curves are the results using folded and $J^+$ decompositions, respectively.}
\label{fig:ex_int_3d}
\end{figure}

Now we integrate Bernstein polynomials over a 3D domain.
The corresponding B-rep is shown in Fig.\ \ref{fig:ex_int_3d}(a), which consists of 4 axis-aligned planar surfaces and a bi-quadratic \bz surface.
The control points of the \bz surface are given as $(1,0.2,0)$, $(1,0.8,0.5)$, $(1,0.4,1)$, $(0.5,0.5,0)$, $(0.5,0.5,0.5)$, $(0.25,0.25,1)$, $(0.2,1,0)$, $(0,1,0.5)$, and $(0.3,1,1)$.
The B-rep is embedded into a $1\times 1\times 1$ Cartesian grid, i.e., a single element.
Similarly, $11\times 11\times 11$ points are sampled as seed vertices to create a family of pyramidal decompositions, where the ones that yield the worst and best integration errors are chosen.
The worst decomposition is guaranteed to have folded cells.
A $J^+$ decomposition is created using the seed vertex as $(1,1,0.5)$.
The degrees of Bernstein polynomials are set to be $p\in\{3,6,9\}$.
The results are summarized in Fig.~\ref{fig:ex_int_3d}(b), and we observe the following.
First, machine precision can be again achieved using both folded and $J^+$ decompositions with a sufficient number of quadrature points per direction.
Second, machine precision is achieved with fewer quadrature points than the ones theoretically predicted for exact integration.
In 3D, according to Eq.~\eqref{eq:n3d} in~\ref{appA}, we need $n\geq 3(p+1)$ to exactly integrate polynomials.
For instance, when $p=12$, 39 quadrature points are needed in theory, but as shown in Fig.~\ref{fig:ex_int_3d}(b), 17 points already suffice to reach machine precision when the $J^+$ decomposition is used.
Third, given the same number of quadrature points, the results using the $J^+$ decomposition fall between those of the best and the worst decompositions. 
Even in the worst scenario, the difference in accuracy can be within two orders provided that a sufficient number of quadrature points are used.
In other words, the performance of a folded decomposition can be always comparable to that of $J^+$ decompositions.

\subsection{Integration of general functions}

In this section, we study to integrate several regular functions that are not polynomials. 
We use again the B-reps shown in Figs.\ \ref{fig:ex_poly_curve}(a) and \ref{fig:ex_int_3d}(a) for 2D and 3D tests, respectively. 
The integrands of interest are 
$f(x,y)=\exp(y)\sin(x)\cos(y)$
and
$f(x,y,z)=\exp(y)\sin(x)\cos(y)\cos(z)$
for 2D and 3D tests, respectively. 
In 2D, the B-rep is embedded into a series of globally refined Cartesian grids, from $2\times 2$ up to $128\times 128$. 
Likewise in 3D, the Cartesian grids range from $2\times 2\times 2$ to $64\times 64\times 64$. 
Integration is performed elementwise. 
The B-reps are approximated in each involved integration cell with quadratic \bz curves in 2D, or bi-quadratic \bz surfaces in 3D.
In each integration cell, we change the number of quadrature points per direction (denoted by $n$) and study the convergence with respect to mesh refinement.
In either 2D or 3D, an overkill solution, which is obtained using a dense Cartesian grid and a large number of quadrature points, serves as the reference solution to evaluate the error.

\begin{figure}[htb]
\centering
\begin{tabular}{cc}
\includegraphics[width=0.47\textwidth]{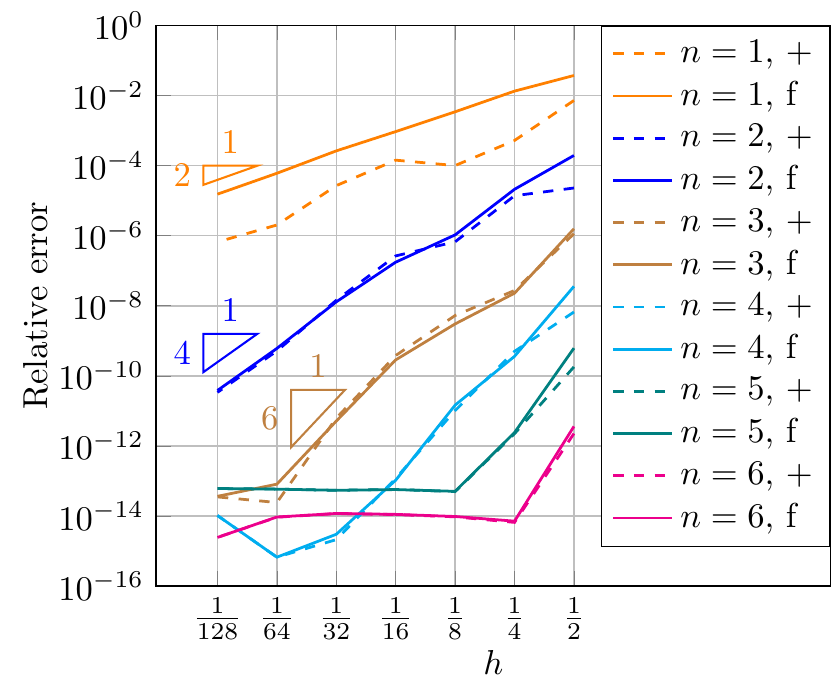} &
\includegraphics[width=0.47\textwidth]{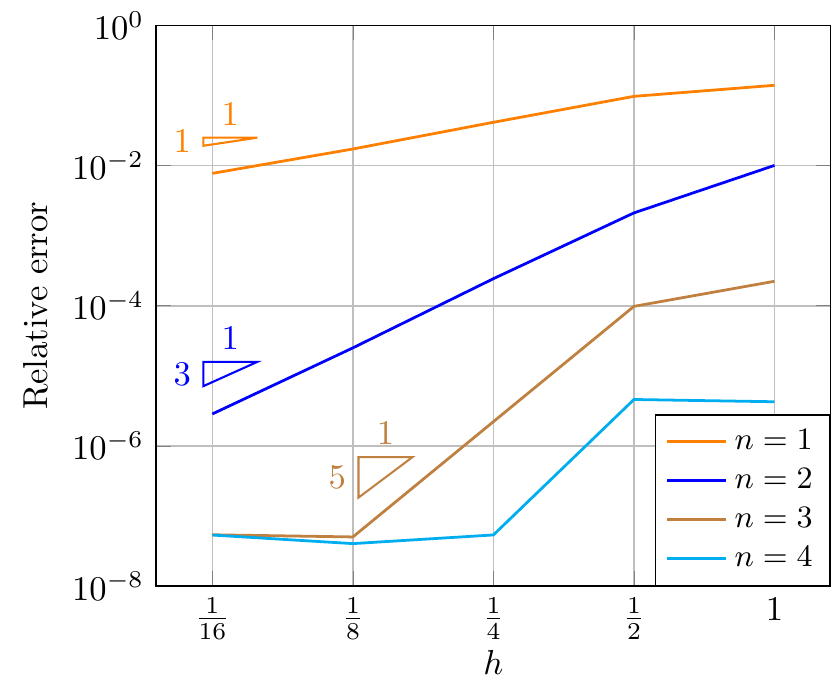} \\
(a) & (b)
\end{tabular}
\caption{Convergence tests of integrating general functions with respect to the mesh size $h$. (a, b) The 2D and 3D tests on the domains shown in Figs.\ \ref{fig:ex_poly_curve}(a) and Fig.\ \ref{fig:ex_int_3d}(a), respectively.
In (a), we compare the performance of folded and $J^+$ decompositions, whereas in (b), we only use folded decompositions.
Solid and dashed curves in (a) are the results using folded and $J^+$ decompositions, respectively.}
\label{fig:conv_int}
\end{figure}

We summarize the error convergence plots in Fig.\ \ref{fig:conv_int}.
We compare folded and $J^+$ decompositions in 2D, whereas we only use folded decompositions in 3D.
In 2D, we observe that the performance of folded and $J^+$ decompositions is comparable, especially when $n\geq 2$.
We also observe a convergence rate of $2n$ when $1\leq n\leq 3$, which is expected using standard Gauss quadrature.
When $n\geq 4$, the integration error quickly goes under machine precision, so we cannot observe the asymptotic convergence of the expected rate.
On the other hand, in 3D (Fig.\ \ref{fig:conv_int}(b)), we observe a convergence rate around $2n-1$ when $1\leq n\leq 3$, one order less than expected.
Recall that a certain surface of an integration cell, which has a \bz representation, approximates the given trimming surface.
The suboptimal convergence may be due to the fact that the geometric error in creating integration cells is more prominent in 3D than that in 2D.

\section{Poisson's problem with immersed IGA}
\label{sec:pde}

In this section, we solve Poisson's problem using IGA, where folded decompositions are used for numerical integration.
We first perform convergence studies on trimmed planar geometries, where the results are compared with those of $J^+$ decompositions.
Moreover, to show the flexibility and potential of folded decompositions, we solve the problem on several trimmed volumes, including an arbitrarily cut cube and a complex engine model.
Note that the number of quadrature points per direction suffices as long as the integration error does not dominate simulation error.
Therefore, we adopt the ($p+r$)-Gauss quadrature rule for a given spline discretization of degree $p$, where we vary $r\geq 1$ to check how it influences convergence.
Moreover, recall that to create integration cells, approximation is often needed with \bz curves (2D) or \bz surfaces (3D).
We choose the degree of each \bz representation to be $p$, which suffices to deliver optimal convergence rates \cite{ref:antolin19a}.

We take Poisson's problem as the model problem in all the following tests. The goal is to find $u\in\Omega$ such that
\begin{subequations}
\begin{alignat}{2}
\Delta u &= f\quad&&\text{in }\Omega\,,\\
\nabla u\cdot\mathbf{n} &= g\quad&&\text{on }\partial\Omega_N\,,\\
u &= 0\quad&&\text{on }\partial\Omega_D\,,
\end{alignat}
\end{subequations}
where $\partial\Omega = \partial\Omega_D \cup \partial\Omega_N$ and $\partial\Omega_D \cap \partial\Omega_N = \emptyset$, $\mathbf{n}$ is the outer unit normal on $\partial\Omega$, and $g$ is the Neumann boundary condition.
We adopt $u(x,y) = \sin(\pi x)\sin(\pi y)$ in 2D and $u(x,y,z) = \sin(\pi x)\sin(\pi y)\sin(\pi z)$ in 3D as manufactured solutions (except for the complex engine model).
The choice of such regular functions as target solutions was motivated by the aim of focusing our study on the consistency error, mainly controled by the numerical integration and geometric represention errors.
The effect of the approximation error in the solution of elliptic PDEs using trimmed spline spaces has been previously addressed in \cite{ref:antolin19a}.

\begin{figure}[htb]
\centering
\begin{tabular}{ccc}
\includegraphics[width=0.33\textwidth]{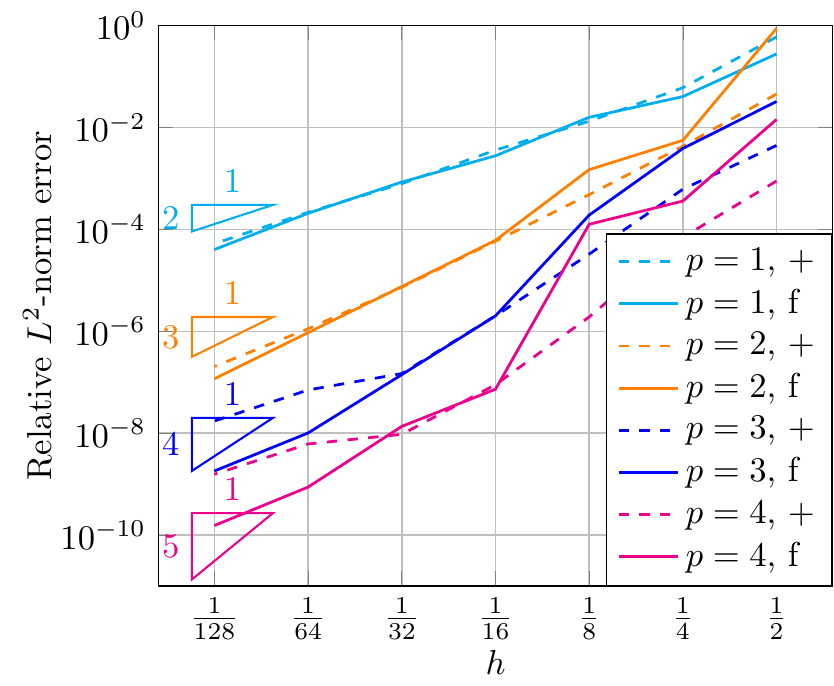} &\hspace{-5mm}
\includegraphics[width=0.33\textwidth]{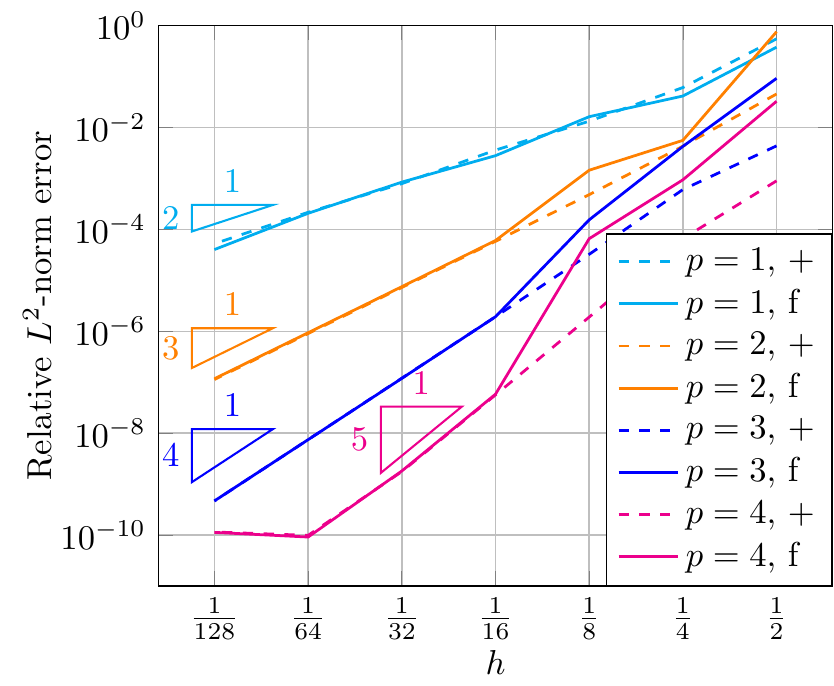} &\hspace{-5mm}
\includegraphics[width=0.33\textwidth]{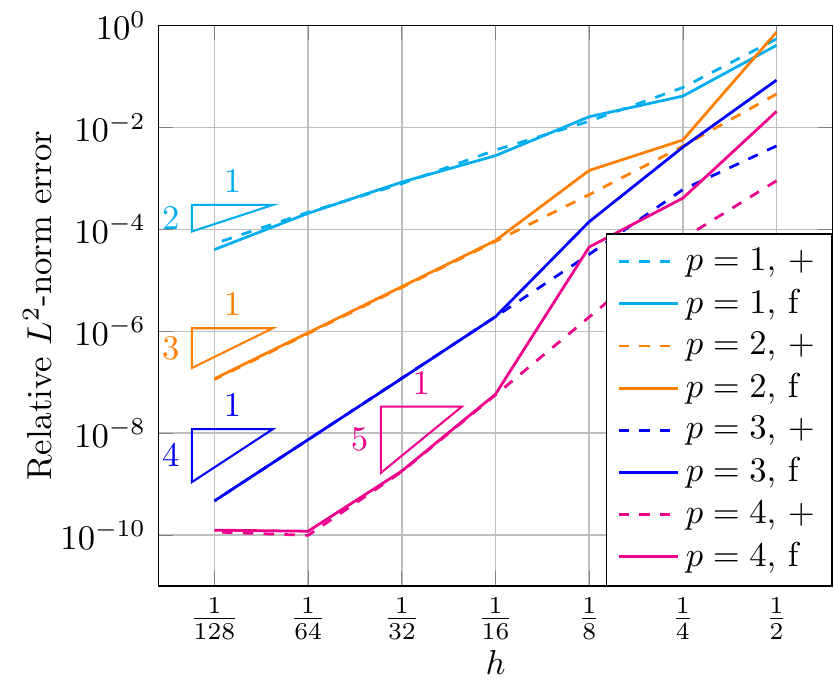} \\
(a) $r=1$ & (b) $r=3$ & (c) $r=5$\\
\includegraphics[width=0.33\textwidth]{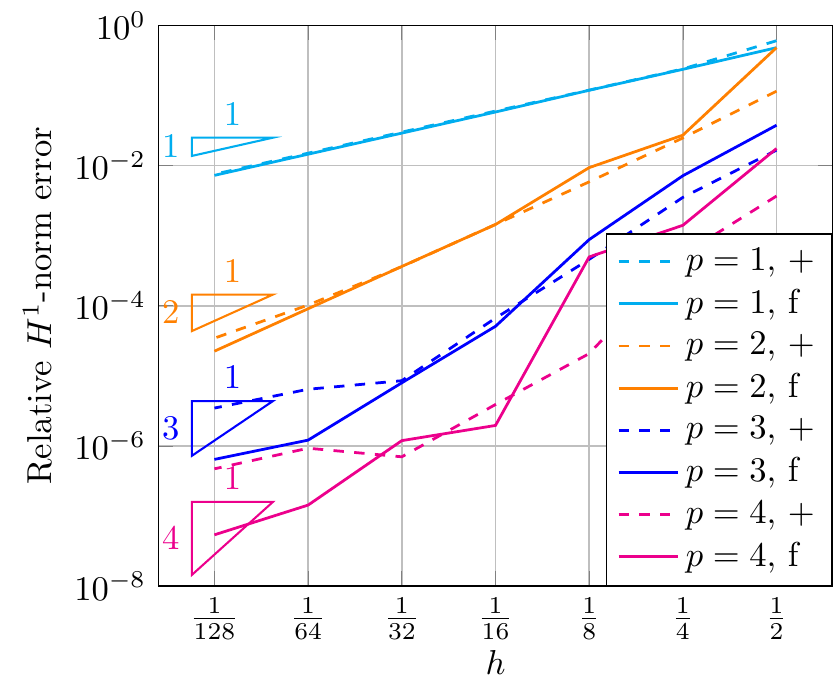} &\hspace{-5mm}
\includegraphics[width=0.33\textwidth]{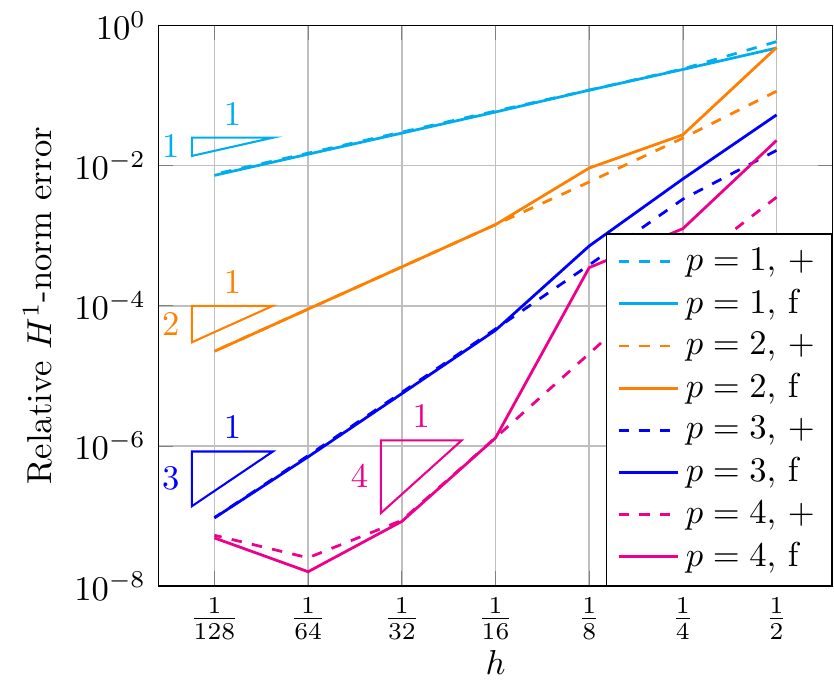} &\hspace{-5mm}
\includegraphics[width=0.33\textwidth]{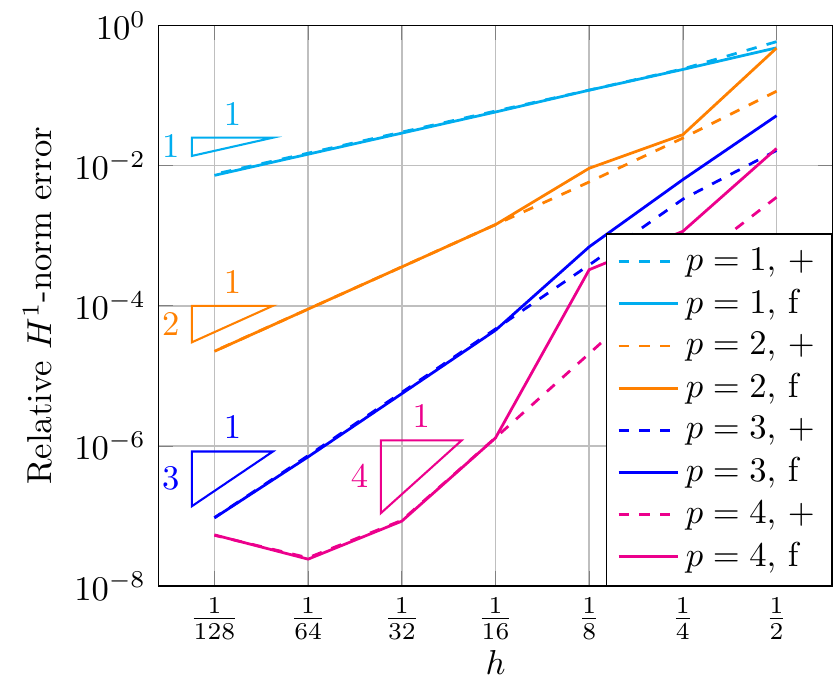} \\
(d) $r=1$ & (e) $r=3$ & (f) $r=5$
\end{tabular}
\caption{The convergence plot in solving Poisson's problem on the domain shown in Fig.\ \ref{fig:ex_poly_curve}(a), where we compare the performance between folded and $J^+$ decompositions for different values of $r$.
(a--c) Relative $L^2$-norm error, and (d--f) relative $H^1$-norm error, where $h$ stands for the element size.
Solid and dashed curves are the results using folded and $J^+$ decompositions, respectively.}
\label{fig:conv_2d_bsp}
\end{figure}

\begin{figure}[htb]
\centering
\begin{tabular}{ccc}
\includegraphics[width=0.33\textwidth]{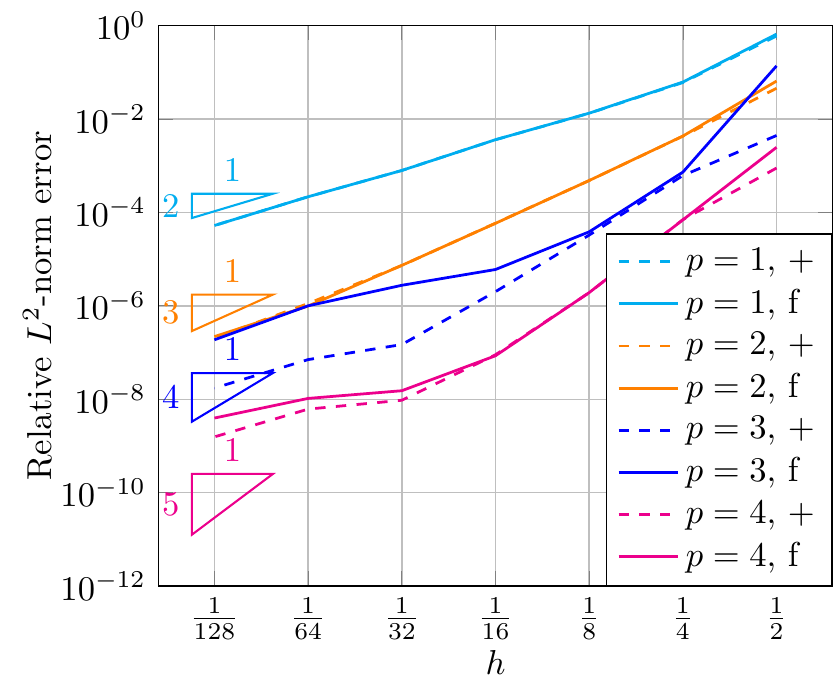} & \hspace{-5mm}
\includegraphics[width=0.33\textwidth]{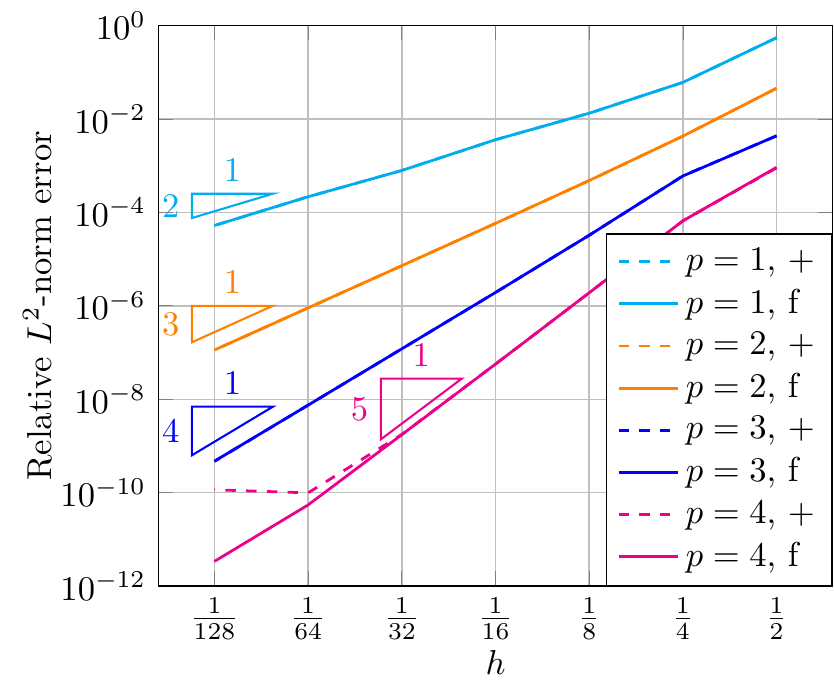} &\hspace{-5mm}
\includegraphics[width=0.33\textwidth]{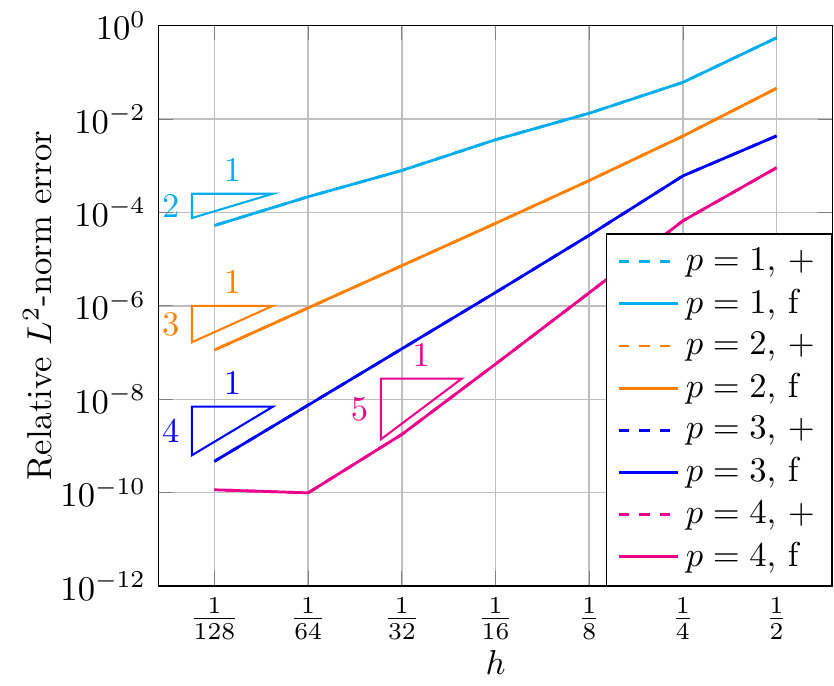} \\
(a) $r=1$ & (b) $r=3$ & (c) $r=5$\\
\includegraphics[width=0.33\textwidth]{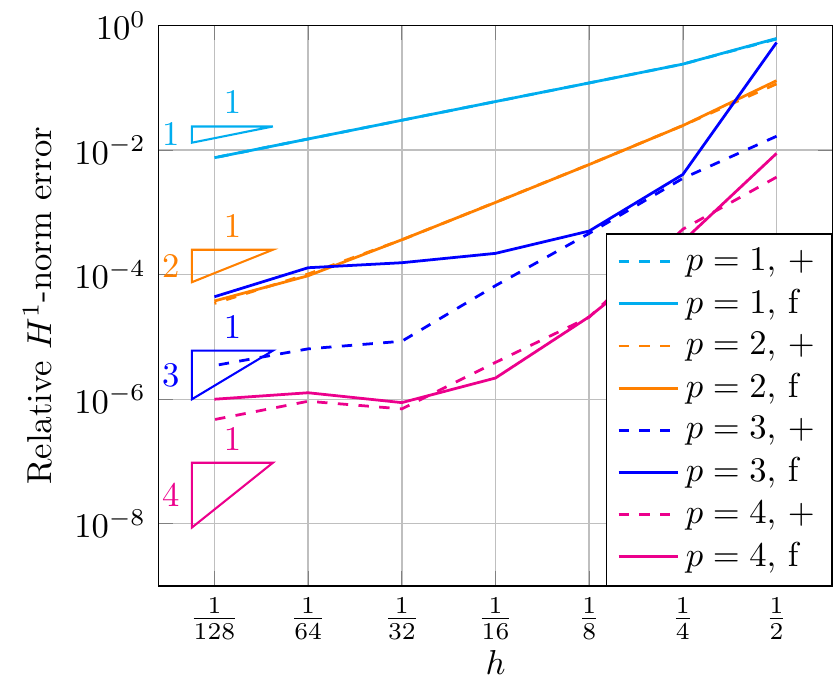} & \hspace{-5mm}
\includegraphics[width=0.33\textwidth]{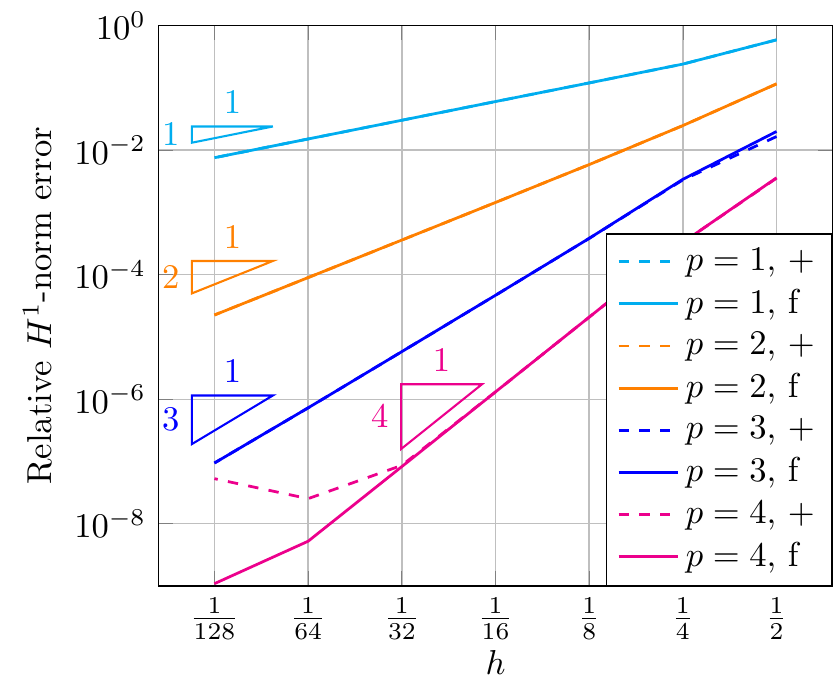} &\hspace{-5mm}
\includegraphics[width=0.33\textwidth]{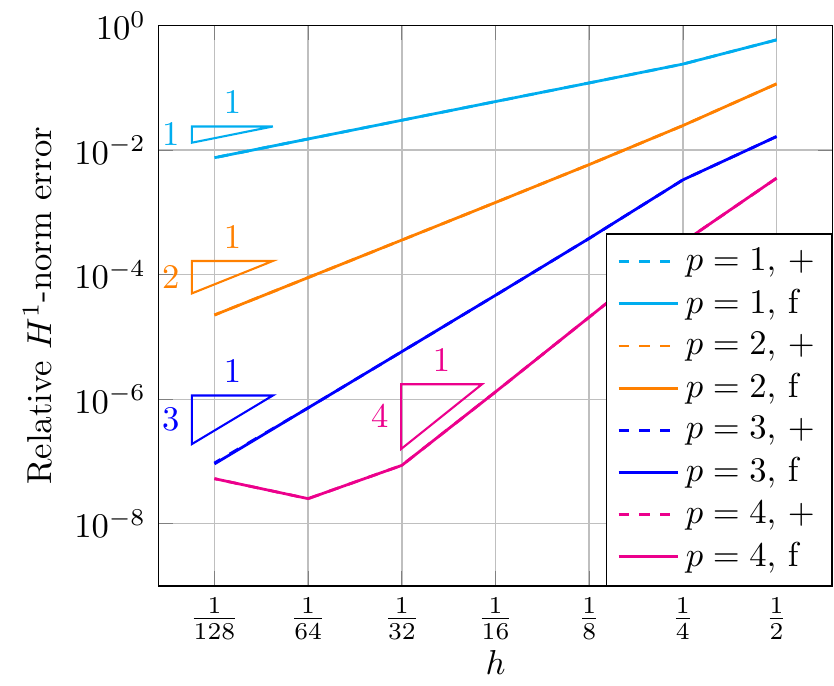} \\
(d) $r=1$ & (e) $r=3$ & (f) $r=5$
\end{tabular}
\caption{The convergence plot in solving Poisson's problem on the domain shown in Fig.\ \ref{fig:ex_ration}(a), where we compare the performance between folded and $J^+$ decompositions for different values of $r$.
(a--c) relative $L^2$-norm error, and (d--f) relative $H^1$-norm error, where $h$ stands for the element size.
Solid and dashed curves are the results using folded and $J^+$ decompositions, respectively.}
\label{fig:conv_2d_hole}
\end{figure}

\subsection{Trimmed planar geometries}

We start with two planar geometries, whose B-reps are already shown in Figs.\ \ref{fig:ex_poly_curve}(a) and \ref{fig:ex_ration}(a).
The problem settings for these two domains are identical.
The homogeneous Dirichlet boundary condition is imposed on the right boundary ($x=1$), whereas the non-homogeneous Neumann boundary condition is imposed elsewhere.
In each problem, the B-rep is embedded into a series of consecutively refined Cartesian grids, from $2\times 2$ up to $128\times 128$.
Different degrees $p$ of B-spline discretizations are defined on each Cartesian grid, with $p\in\{1,2,3,4\}$.
In each integration cell, $p+r$ quadrature points are used along each direction, where we use $r\in\{1,3,5\}$.

For these two geometries, the folded decomposition for each trimmed element is created using the element's origin as seed vertex for triangulation.
This is meant to create highly distorted cells and it guarantees the presence of folded cells for this particular geometry.
$J^+$ decompositions are generated using the quadrangulation method proposed in~\cite{ref:wei21b}.

The convergence plots are shown in Figs.\ \ref{fig:conv_2d_bsp} and \ref{fig:conv_2d_hole}.
We observe the following.
First, with a sufficient number of quadrature points (e.g., $r=5$), we can always achieve optimal convergence rates using both folded and $J^+$ decompositions.
On the other hand, convergence is deteriorated when quadrature points are not enough (e.g., $r=1$).
Second, the convergence curves with folded and $J^+$ decompositions eventually become identical when enough quadrature points are used; see, for instance, when $h\leq 1/16$ and $r=5$.
Third, different values of $r$ are needed for different spline discretizations to achieve optimal convergence.
For instance, with the bi-linear and bi-quadratic spline discretizations ($p=1$ and $p=2$), $r=1$ suffices to retain optimal convergence.
However, when $p=4$, we need $r\geq 3$.
We also observe in this case that the convergence curve hits a plateau when $H^1$-norm error reaches $10^{-8}$.
This is caused by the geometric tolerances involved in the computation of surface-surface intersections required to locate trimmed elements, as also discussed in~\cite{ref:antolin21}.
Once the simulation result reaches this level of accuracy, the geometric approximation error starts to dominate and thus leads to a plateau in convergence.
These tests provide the numerical evidence that as long as a sufficient number of quadrature points are used, folded decompositions can yield the same convergence as $J^+$ decompositions.

\subsection{Trimmed volumes}

\begin{figure}[!htb]
\centering
\begin{tabular}{ccc}
\includegraphics[width=0.28\textwidth]{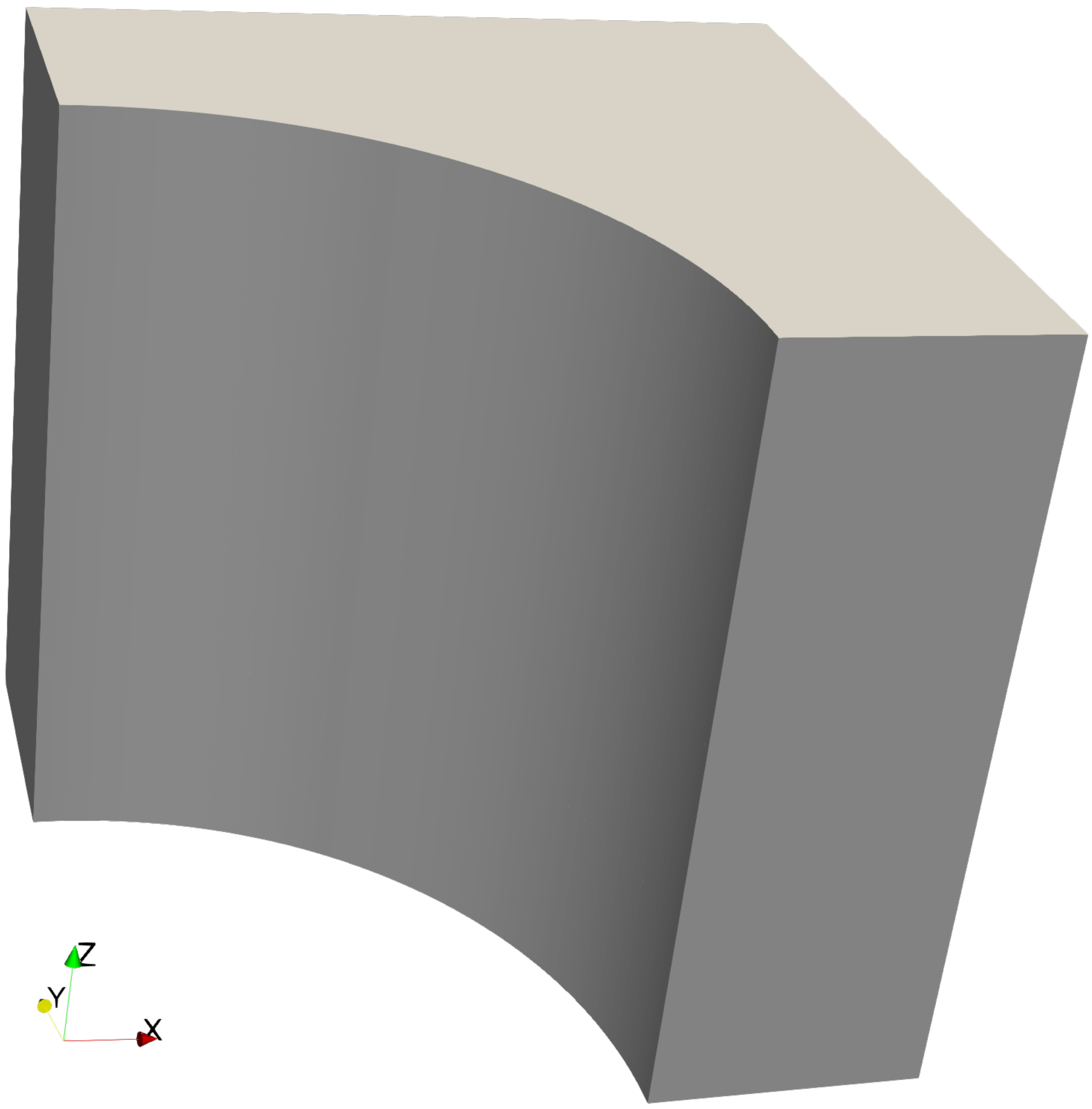} &\hspace{-5mm}
\includegraphics[width=0.35\textwidth]{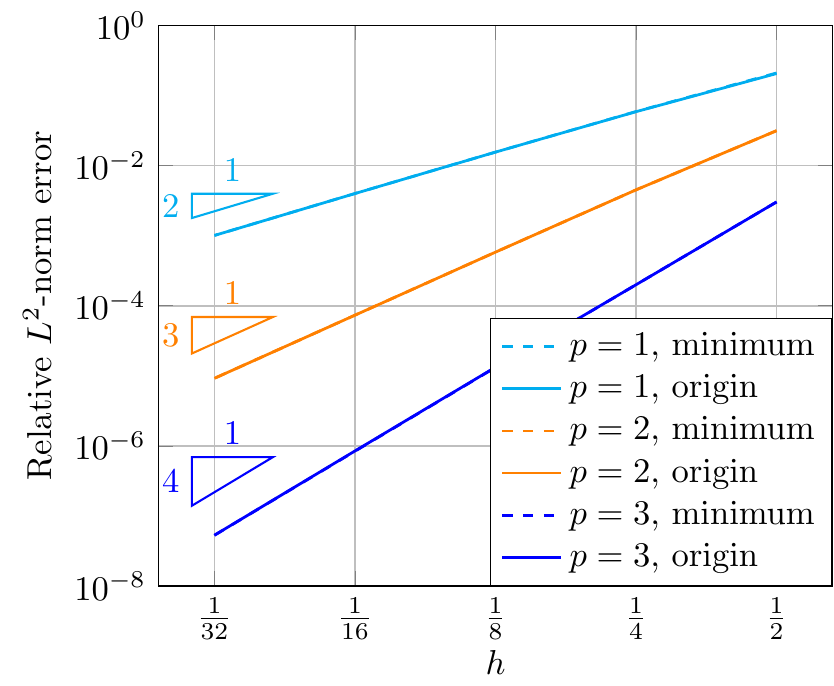} &\hspace{-5mm}
\includegraphics[width=0.35\textwidth]{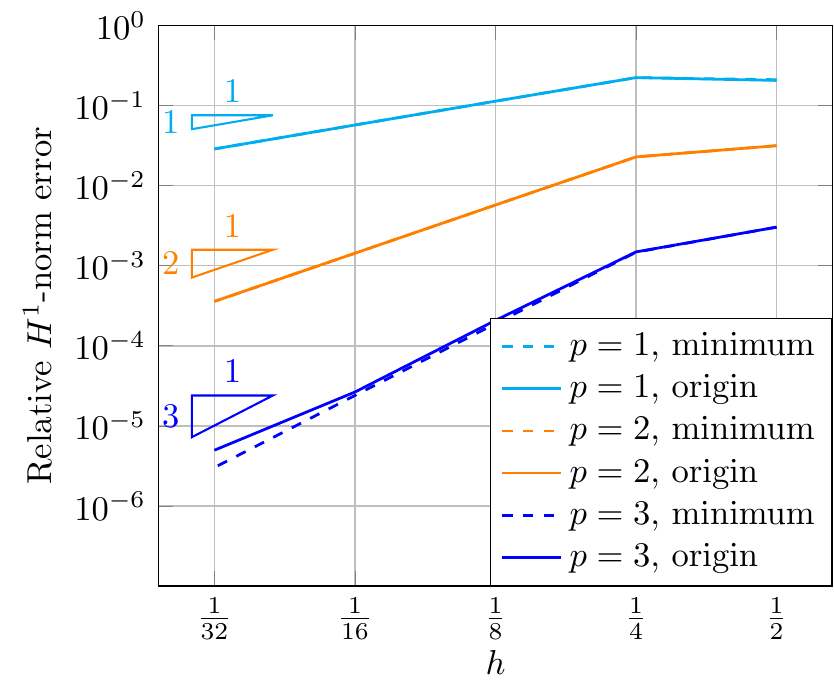}\\
(a) & (b) & (c)
\end{tabular}
\caption{Solving Poisson's problem on a domain obtained by cutting a cube with a cylinder.
(a) The input B-rep, and (b, c) convergence plots of relative $L^2$- and $H^1$-norm errors, respectively, using $r=5$ and two different ways of decompositions: one leading to the minimum number of cells and the other choosing the element origin as the seed vertex.}
\label{fig:conv_3d_cylin}
\end{figure}

Now we present several 3D examples where we do not constrain the type of decompositions to be generated. 
In other words, the decomposition algorithm terminates as long as a decomposition is achieved for a trimmed element, which can be either a $J^+$ decomposition or a folded decomposition. 
Nonetheless, folded decompositions are often generated in practice.

We start with a simple geometry as shown in Fig.~\ref{fig:conv_3d_cylin}(a), where a cube $[0,L]^3$ ($L=1000$) is cut by a cylinder of radius $0.65L$, and the axis of the cylinder is along the $z$ direction.
A quarter of the cylinder is involved in the B-rep, which is exactly represented by a rational bi-quadratic \bz surface.
To perform a convergence study, the B-rep is embedded into a series of consecutively refined Cartesian grids, with the resolution from $2\times 2\times 2$ up to $32\times 32\times 32$.
The homogeneous Dirichlet boundary is imposed on the (bottom) plane $x=L$, whereas the non-homogeneous Neumann boundary condition is imposed elsewhere.
On each Cartesian grid, we define B-spline discretizations of degrees $p\in\{1,2,3\}$.
When creating integration cells for each trimmed element, we consider two options for the seed vertex: (1) one of the polyhedral vertices (recall that a trimmed element is a curved polyhedron), and (2) the origin of each element.
Among all the possible decompositions in Option (1), we choose the one that leads to the minimum number of cells.
Option (2) is meant to create highly distorted cells and it guarantees the presence of folded cells for this particular geometry.
In each cell, $p+5$ Gauss quadrature points are used per direction.
The convergence results are summarized in Figs.~\ref{fig:conv_3d_cylin}(b, c), where we compare the two options of decompositions.
We observe that they yield almost indistinguishable convergence results.
Optimal convergence rates are achieved with all the B-spline discretizations.
This implies that the accuracy and convergence are independent of how we create decompositions.


\begin{figure}[htb]
\centering
\begin{tabular}{cc}
\includegraphics[height=0.35\textwidth]{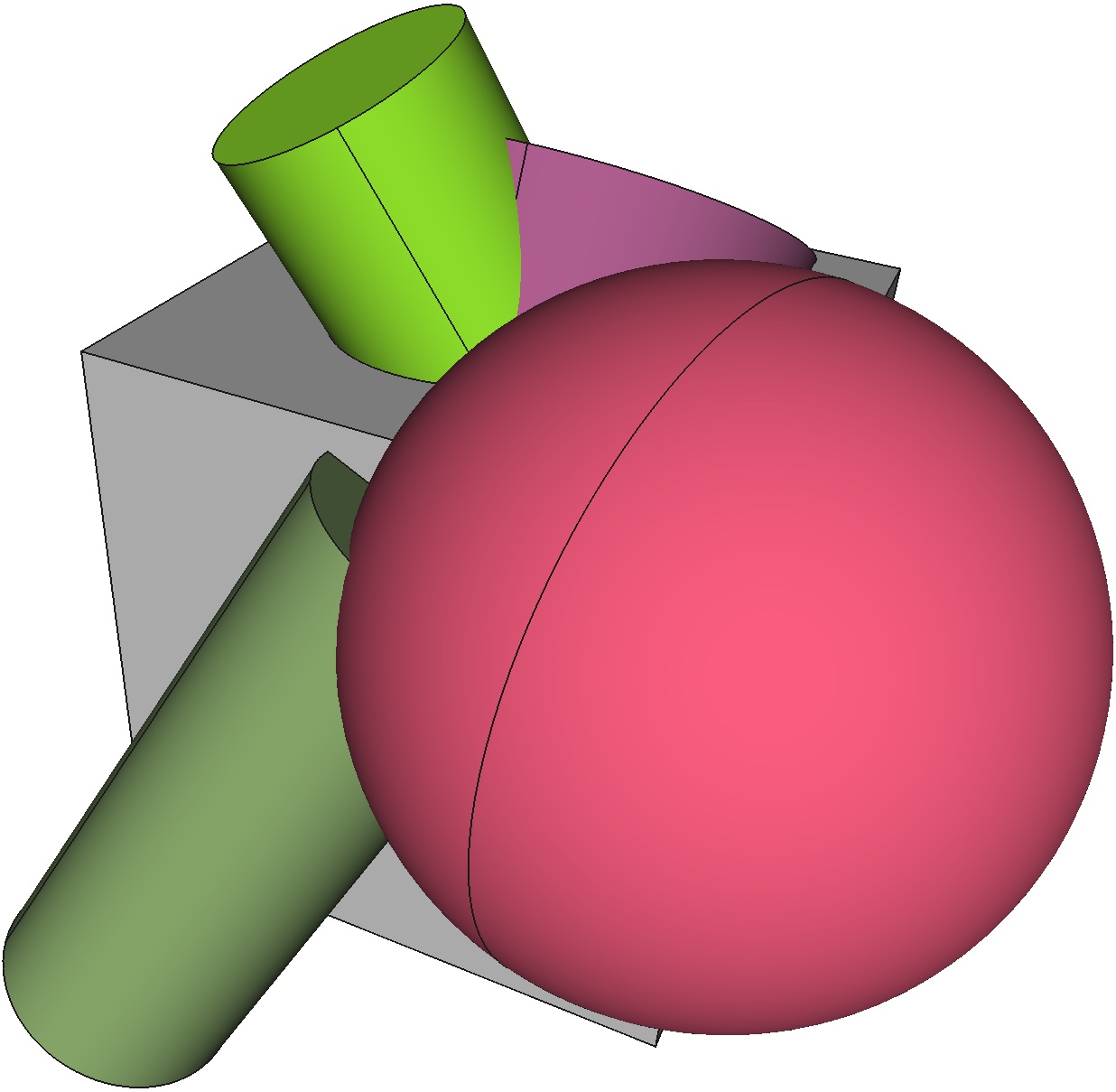} &\hspace{+2mm}
\includegraphics[height=0.35\textwidth]{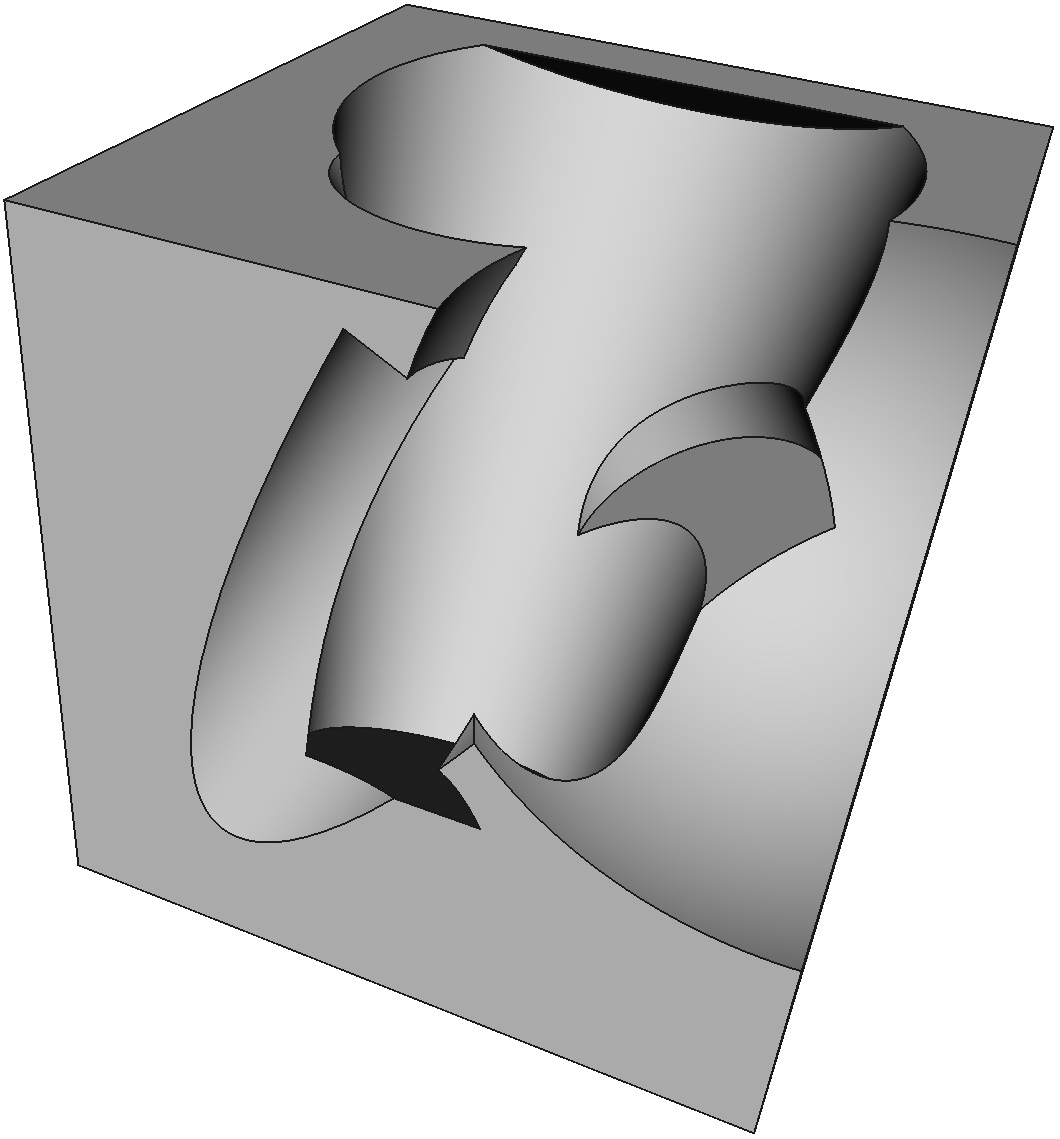}\\
(a) & (b)
\end{tabular}
\caption{Construction of a complex geometry through trimming.
(a) A cube to be cut by two cylinders (green), a cone (violet), and a sphere (pink), and (b) the resulting B-rep.}
\label{fig:complex_geom_inp}
\end{figure}

\begin{figure}[!htb]
\centering
\begin{tabular}{ccc}
\includegraphics[height=0.33\textwidth]{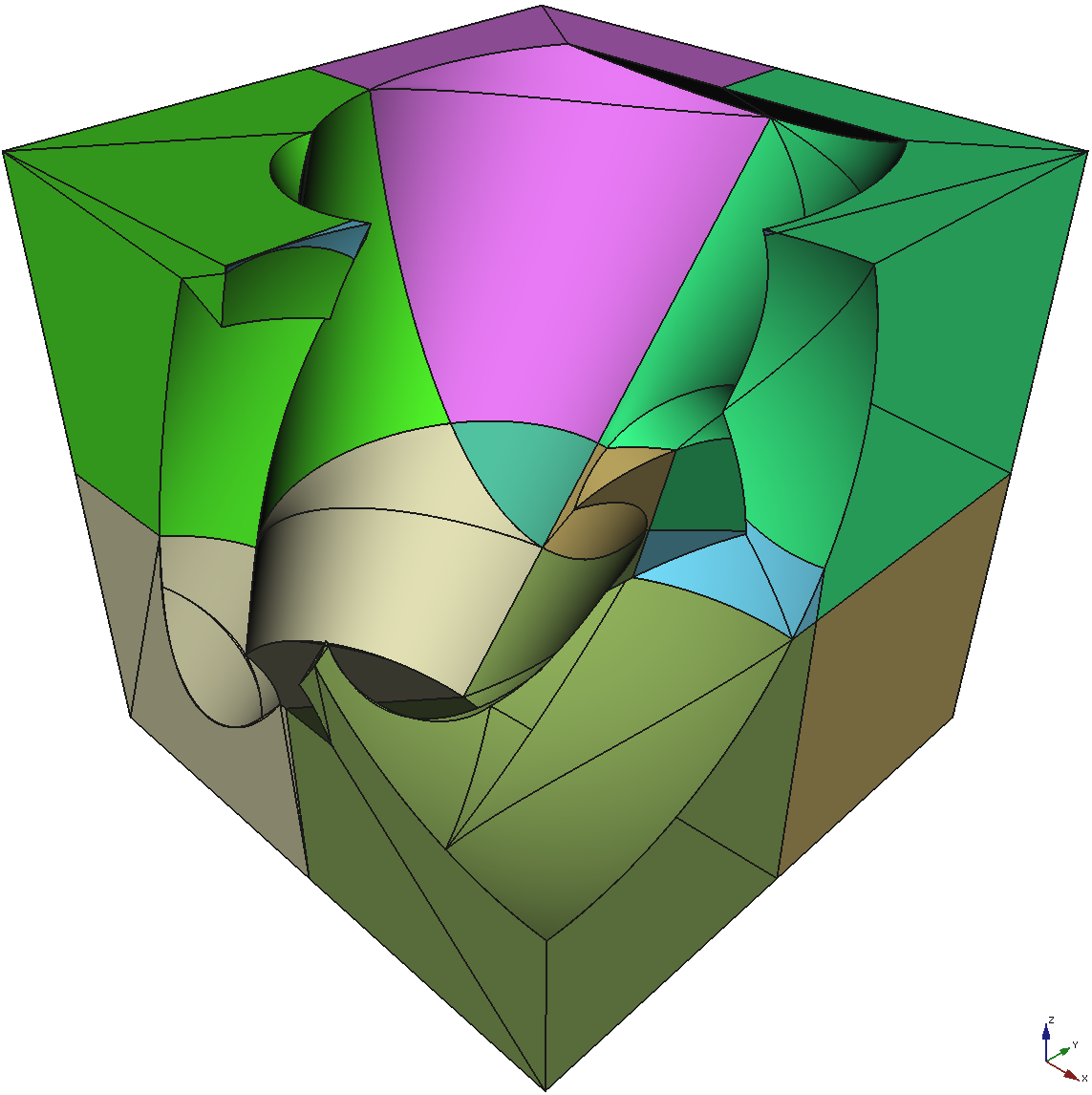} &\hspace{-2mm}
\includegraphics[height=0.33\textwidth]{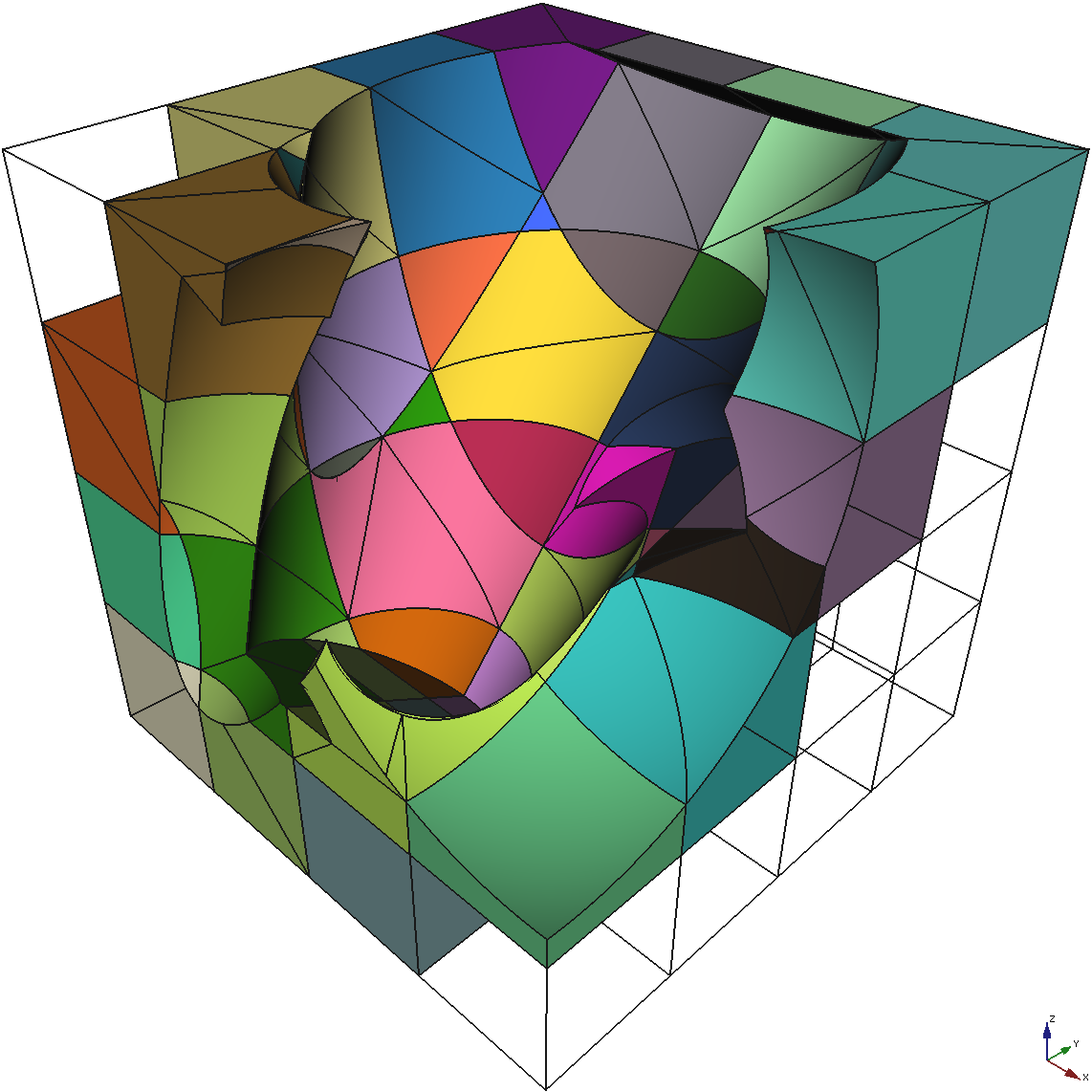}&\hspace{-2mm}
\includegraphics[height=0.33\textwidth]{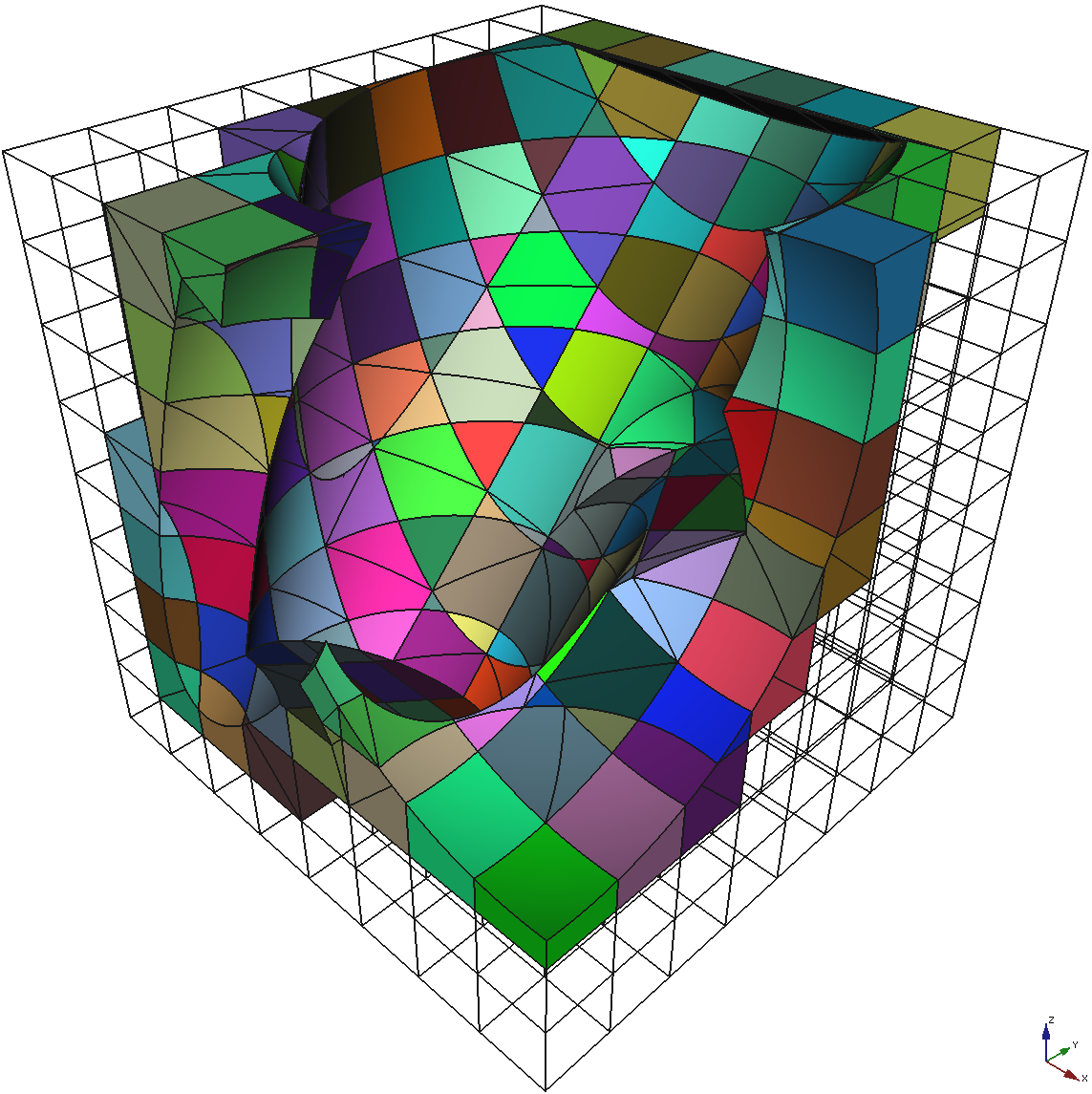}\\
(a) $2\times 2\times 2$ & (b) $4\times 4\times 4$ & (c) $8\times 8\times 8$ \\
\end{tabular}
\begin{tabular}{cc}
\includegraphics[height=0.33\textwidth]{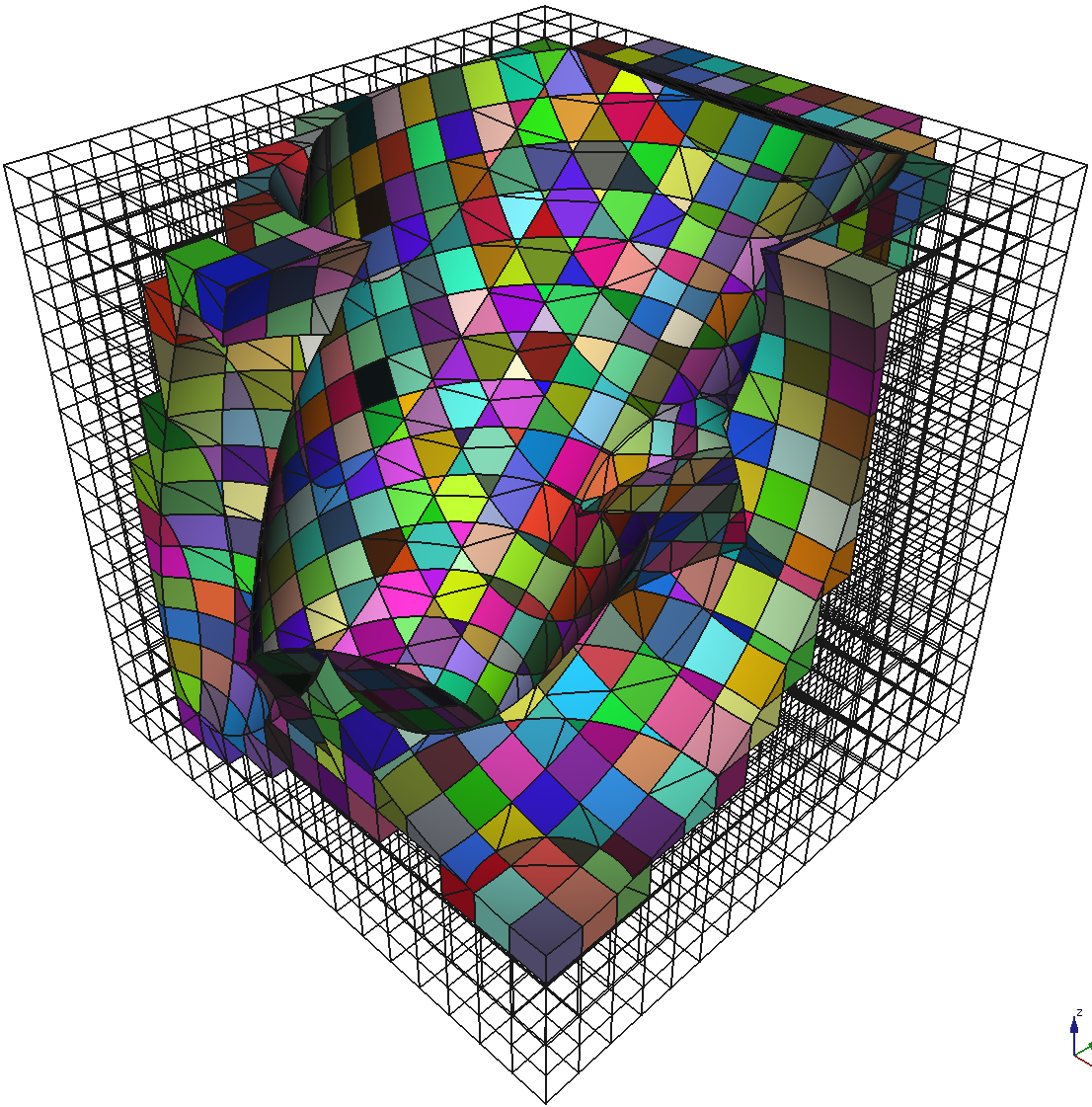}&\hspace{-2mm}
\includegraphics[height=0.33\textwidth]{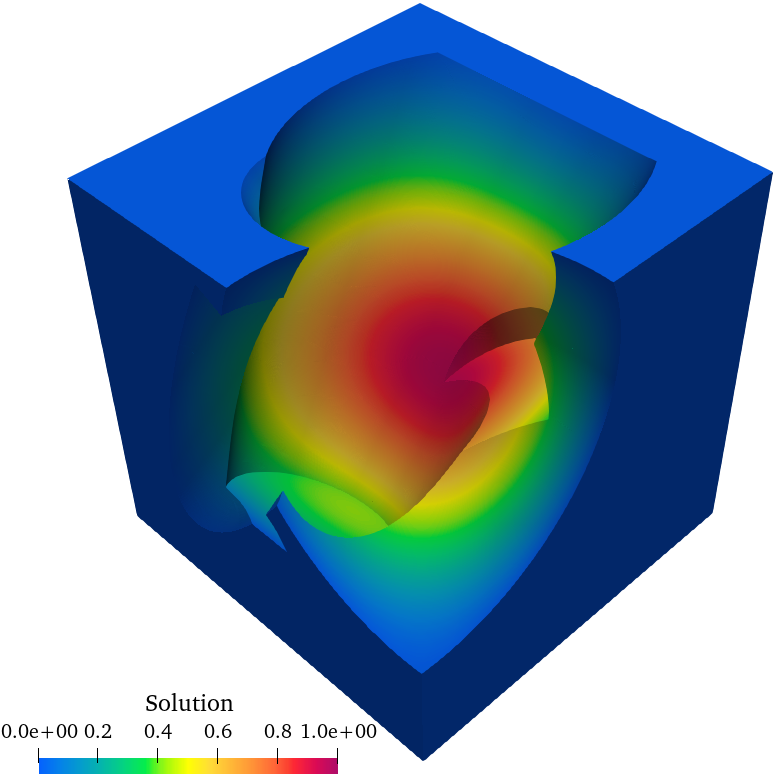}\\
(d) $16\times 16\times 16$ & (e) Solution on (d)
\end{tabular}
\caption{The B-rep embedded into a series of consecutively refined Cartesian grids (a--d), and the solution field (e) obtained using the mesh in (d).
In (a--d), only trimmed elements are shown and each color represents a trimmed element.
The black curves represent face-wise decompositions for the boundaries of trimmed elements as well as the outlines of non-trimmed elements.}
\label{fig:complex_geom_mesh}
\end{figure}

Next, we consider a complex geometry shown in Fig.~\ref{fig:complex_geom_inp}, where a cube $[0,L]^3$ ($L=10$) is cut by a sphere, a cone, and two cylinders.
The corresponding B-rep is shown in Fig.~\ref{fig:complex_geom_inp}(b), which contains several sharp curves/corners.
Such features can be respected with integration cells for accurate integration.
When creating a decomposition for each trimmed element, the seed vertex is chosen to be one of the element vertices.
Among all the possible decompositions, we choose the one that leads to the minimum number of cells.
The homogeneous Dirichlet boundary is imposed on the plane $z=0$, whereas the non-homogeneous Neumann boundary condition is imposed elsewhere.
All the other problem settings, including Cartesian grids, corresponding spline discretizations, and the number of quadrature points per direction, follow the same as those in the previous example.

\begin{figure}[!htb]
\centering
\begin{tabular}{cc}
\includegraphics[width=0.35\textwidth]{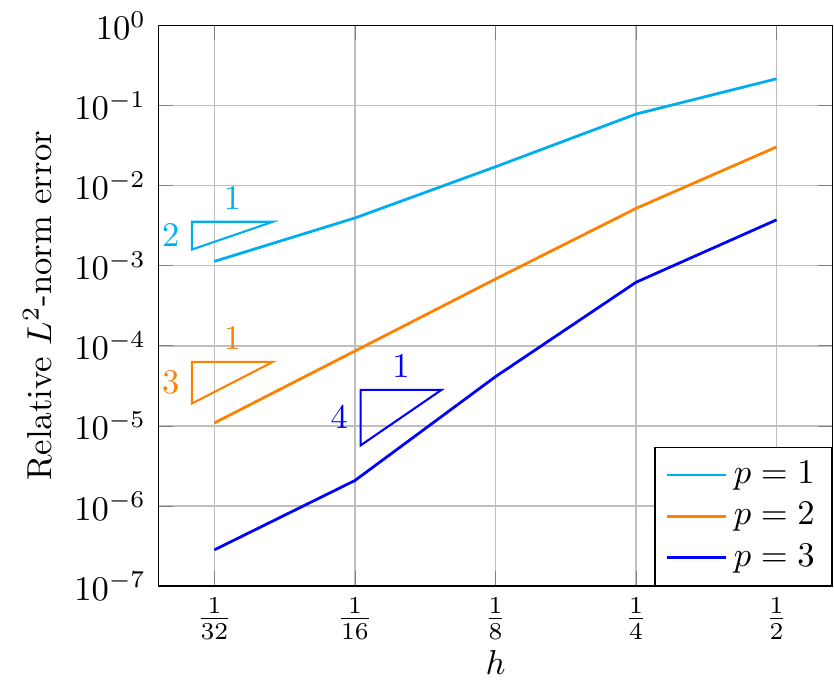} &\hspace{+2mm}
\includegraphics[width=0.35\textwidth]{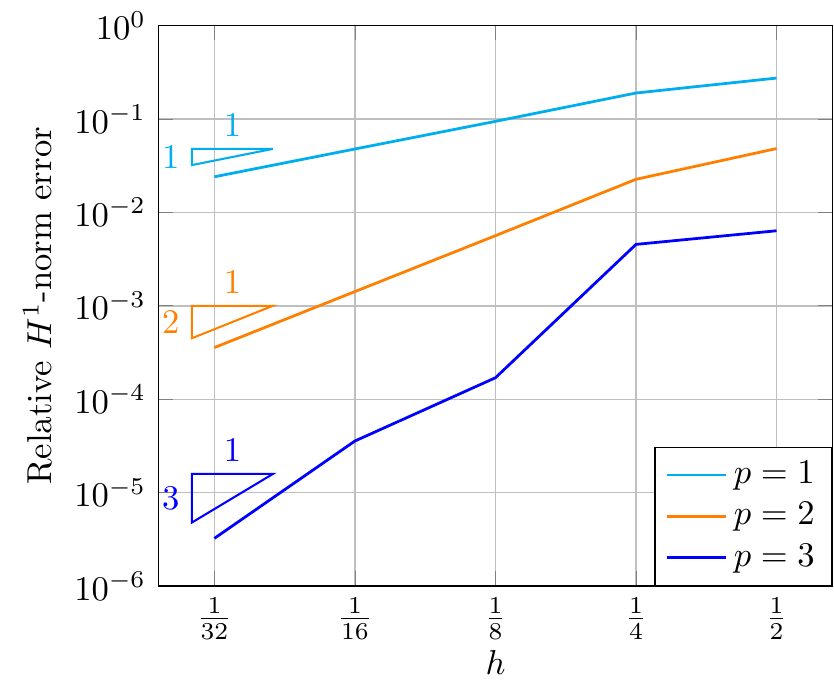}\\
(a) & (b)
\end{tabular}
\caption{Solving Poisson's problem on a domain shown in Fig.\ \ref{fig:complex_geom_inp}.
(a, b) Convergence plots of relative $L^2$- and $H^1$-norm errors, respectively, with $r=5$.}
\label{fig:conv_3d_comp}
\end{figure}

In Fig.~\ref{fig:complex_geom_mesh}, we show a series of Cartesian grids and the corresponding trimmed elements.
Each color represents a single trimmed element, which is generally decomposed into a collection of reparameterized cells for integration.
We observe that such cells well align with the sharp features in the input B-rep.
In addition, in Figs.~\ref{fig:complex_geom_mesh}(a-d) it can be seen that face-wise decompositions are not conformal.
The solution field using the $16\times 16\times 16$ Cartesian grid is shown in Fig.~\ref{fig:complex_geom_mesh}(e).
Moreover, the convergence plots are shown in Fig.~\ref{fig:conv_3d_comp}, where we again observe optimal convergence rates for linear, quadratic and cubic spline discretizations.

\begin{figure}[!htb]
\centering
\begin{tabular}{c}
\includegraphics[width=\textwidth]{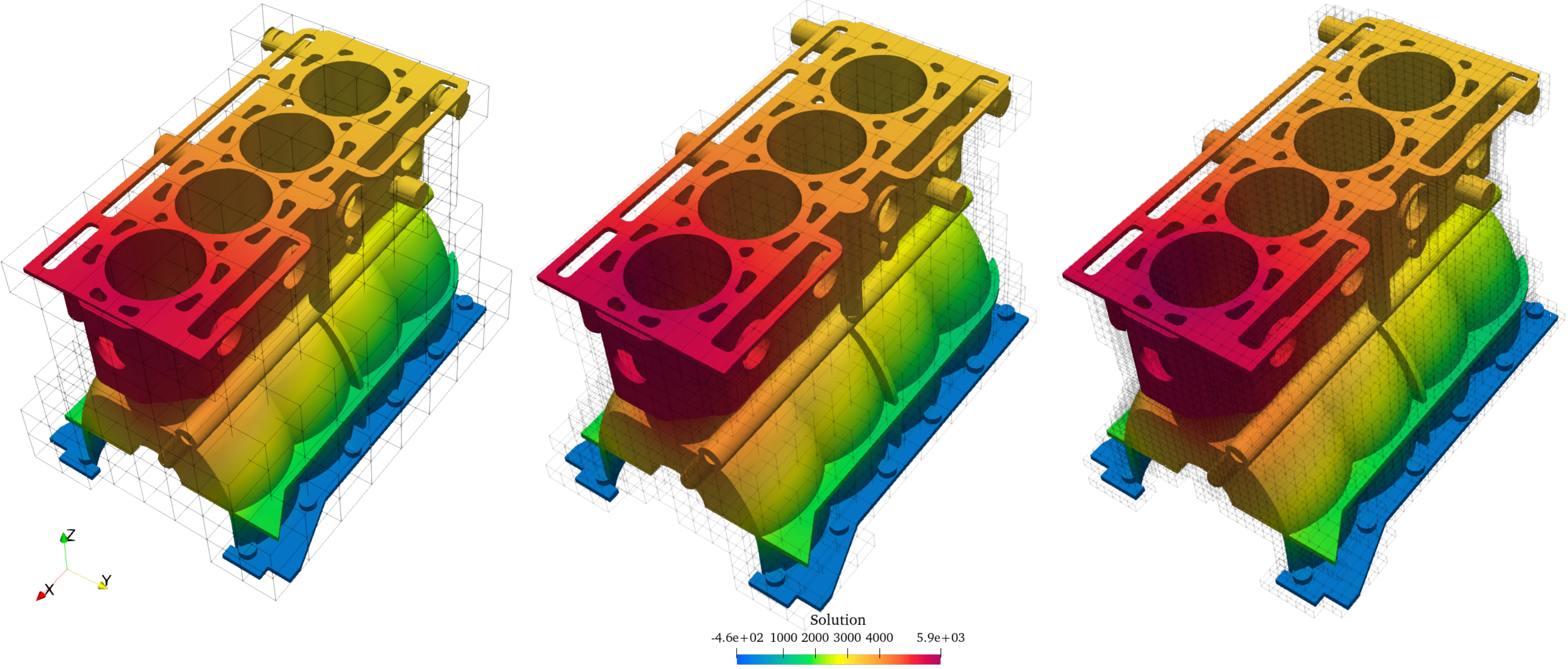}\\
\end{tabular}
\begin{tabular}{ccc}
(a) $8\times 8\times 8$ & \hspace{+0.2\textwidth} (b) $16\times 16\times 16$ & \hspace{+0.2\textwidth} (c) $32\times 32\times 32$
\end{tabular}
\begin{tabular}{c}
\includegraphics[width=\textwidth]{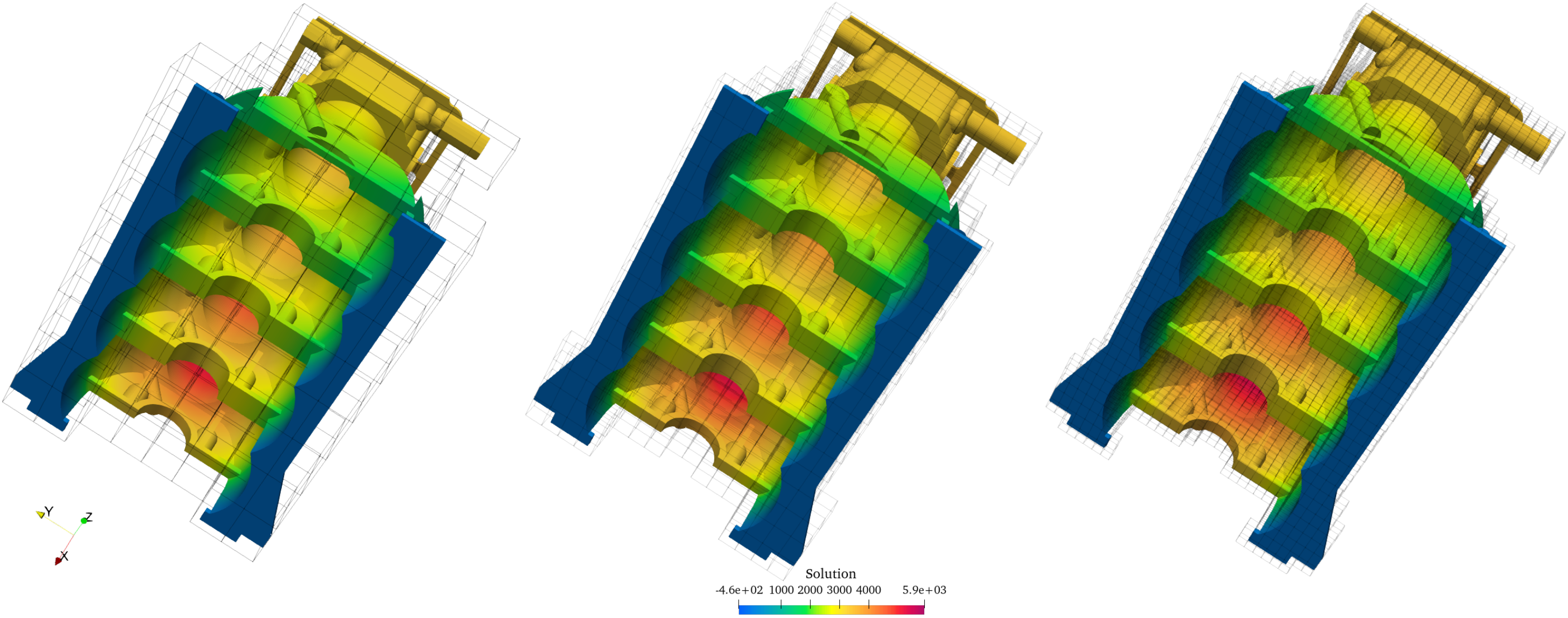}\\
\end{tabular}
\begin{tabular}{ccc}
(d) $8\times 8\times 8$ & \hspace{+0.2\textwidth} (e) $16\times 16\times 16$ & \hspace{+0.2\textwidth} (f) $32\times 32\times 32$
\end{tabular}
\caption{Solving Poisson's problem on an engine model with different resolutions of Cartesian grids.
(a--b) show a top view of the solution for different discretizations, whereas (d--f) show the solution seen from the bottom.
The underlying embedding domain is represented with black lines.}
\label{fig:conv_3d_engine}
\end{figure}

As the last example, we present an engine model obtained from the gallery of OpenCASCADE \cite{ref:opencascade} to show that our proposed method can robustly accommodate real-world geometries.
Its B-rep has many detailed features that are represented by 918 trimmed spline surfaces in total.
The B-rep is embedded into a series of Cartesian grids with resolutions $8\times 8\times 8$, $16\times 16\times 16$, and $32\times 32\times 32$, on which we solve Poisson's problem using tricubic B-splines.
The numerical solutions are shown in Fig.\ \ref{fig:conv_3d_engine}.
This demonstrates that with folded decompositions, we can flexibly reparameterize trimmed elements while retaining all the important geometric features.

\section{Conclusions and future work}
\label{sec:con}

In this paper, we propose to use folded decompositions for numerical integration over complex domains, which essentially eliminates constraints upon candidate decompositions and thus provides a significant flexibility to accommodate real-world geometries in engineering simulation, while retaining a relatively small number of integration cells.
A domain of interest usually has a B-rep in 3D.
A corresponding volumetric decomposition is created following two steps.
First, the B-rep is decomposed into smaller surface patches, which respect possible sharp geometric features and do not need to be conformal.
Such surface patches are generally approximated with \bz surfaces.
Second, the domain is decomposed into a collection of pyramidal cells using the surface patches and a seed vertex chosen in a rather arbitrary manner.
The resulting pyramids serve as the integration cells where standard Gauss quadrature can be applied.
Generally, such cells present negative Jacobian because there is no constraint on the choice of the seed vertex.

The validity of these unconstrained decompositions is theoretically guaranteed.
We have performed two groups of numerical tests to numerically assess the validity and robustness of folded decompositions: integration of functions and solving partial differential equations (PDEs) using isogeometric analysis (IGA).
In both cases, we compare the performance of folded decompositions with that of $J^+$ decompositions.
In the case of integrating polynomials, we have shown that with a sufficient number of quadrature points, machine precision can always be reached using folded decompositions.
On the other hand, in the case of solving PDEs, we observe that the performance of folded and $J^+$ decompositions is almost indistinguishable in terms of accuracy and convergence, as long as a sufficient number of quadrature points are used.
However, due to its great flexibility, folded decompositions can accommodate very complex geometries that involve various kinds of geometric features.

In summary, the proposed methodology allows the creation of high-order decompositions of B-reps for integration purposes in a very simple way, while yielding results analogous in terms of accuracy to the ones obtained using $J^+$ decompositions.
Nevertheless, while folded decompositions provide a significant flexibility and robustness to accommodate real-world complex geometries, there are limitations, for instance, when we use them to visualize simulation results.
Therefore, in the future, we plan to further exploit possible decomposition schemes such that we can minimize the number of integration cells while being able to deliver $J^+$ cells at its best.

\section*{Acknowledgements}

P.\ Antolin and A.\ Buffa are partially supported by the ERC AdG project CHANGE n.\ 694515 and the SNSF through the project ``Design-through-Analysis (of PDEs): the litmus test” n.\ 40B2-0 187094 (BRIDGE Discovery 2019).
X.\ Wei is partially supported by the ERC AdG project CHANGE n.\ 694515 and the SNSF project HOGAEMS n.\ 200021-188589.


\appendix

\section{Number of quadrature points per direction} \label{appA}
In this Appendix we study the number of quadrature points needed per direction to \emph{exactly} integrate a given polynomial $f(\mathbf{x})$ over a certain domain $\Omega$, i.e., to compute $I=\int_{\Omega} f(\mathbf{x})$ (Eq.\ \eqref{eq:integral}).
Let $\Q_{p,p}(\Omega)$ and $\Q_{p,p,p}(\Omega)$ denote bi-degree-$p$ and tri-degree-$p$ polynomial spaces defined on $\Omega$, respectively, and we assume $f(\mathbf{x})\in\Q_{p,p}(\Omega)$ in 2D or $f(\mathbf{x})\in\Q_{p,p,p}(\Omega)$ in 3D.

We adopt the following three assumptions to facilitate the study of exact polynomial integration, which, however, are not needed elsewhere in the paper.
First, we assume that $\Omega$ can be exactly represented by a set of cells, $\fundef{\Tb_\ell}{\Lambda=[0,1]^d}{\Tc_\ell}$.
Second, $\Tb_\ell$ is assumed to be a \bz patch that has a polynomial representation.
In other words, $\mathbf{C}_i(u)$ in \eqref{eq:cell2d} and $\mathbf{S}_i\circ\mathbf{R}_{i,j}(u,v)$ in \eqref{eq:cell3d} are \bz curves and surfaces, respectively.
Third, for simplicity, $\Tb_\ell$ is further assumed to have an identical degree ($q$) in all directions, i.e., $\Tb_\ell\in \Q_{q,q}(\Lambda)\times \Q_{q,q}(\Lambda)$ for $d=2$ or $\Tb_\ell\in \Q_{q,q,q}(\Lambda)\times \Q_{q,q,q}(\Lambda)\times \Q_{q,q,q}(\Lambda)$ for $d=3$.

Therefore, according to~\eqref{eq:decompquad}, we have
\begin{equation}
I=\int_{\Omega} f(\xb) = \sum_{\ell=1}^{L} \int_{\Lambda} f \circ \Tb_\ell \, \jacob{\Tb_\ell}.
\end{equation}
When $d=2$, $f \circ \Tb_\ell\in\Q_{2pq,2pq}(\Lambda)$,  $\jacob{\Tb_\ell}\in\Q_{2q-1,2q-1}(\Lambda)$, and thus we have 
\begin{equation}
f \circ \Tb_\ell \jacob{\Tb_\ell} \in\Q_{2q(p+1)-1,2q(p+1)-1}(\Lambda).
\label{eq:fcomp}
\end{equation}
To exactly compute \eqref{eq:fcomp} with the tensor-product Gauss quadrature defined in $\Lambda$, we need the number of quadrature points per direction $n$ to be sufficiently large such that
\begin{equation}
2n-1 \geq 2q(p+1)-1 \iff n \geq q(p+1).
\label{eq:n2d}
\end{equation}
Likewise when $d=3$, $f \circ \Tb_\ell \in\Q_{3pq,3pq,3pq}(\Lambda)$ and $\jacob{\Tb_\ell} \in\Q_{3q-1,3q-1,3q-1}(\Lambda)$.
We have
\begin{equation}
f \circ \Tb_\ell \jacob{\Tb_\ell} \in\Q_{3q(p+1)-1,3q(p+1)-1,3q(p+1)-1}(\Lambda) ,
\end{equation}
and
\begin{equation}
2n-1 \geq 3q(p+1)-1 \iff n \geq \frac{3}{2} q(p+1).
\label{eq:n3d}
\end{equation}

Eqs.\ (\ref{eq:n2d}, \ref{eq:n3d}) imply that for an exact integration, $n$ can easily become very large when either $p$ or $q$ increases, due to the composition $f\circ \Tb_\ell$.
However, exact integration of a polynomial is often not needed as long as the integration error does not dominate the application of interest.
We have such examples in Section \ref{sec:pde}.

\bibliographystyle{unsrt}
\bibliography{ref}
\end{document}